\documentclass[a4paper,12pt]{article}

\usepackage{times, url}
\usepackage{amsmath,amssymb,amsthm,amsfonts,bm,float,cite}
\usepackage{mathtools}
\usepackage{boldline}
\usepackage[utf8]{inputenc}
\IfFileExists{esvect.sty}{\usepackage{esvect}}{}

\newtheorem{definition}{Definition}[section]

\usepackage[none]{hyphenat}
\usepackage[center]{caption}
\DeclarePairedDelimiter\ceil{\lceil}{\rceil}
\DeclarePairedDelimiter\floor{\lfloor}{\rfloor}
\usepackage[pdftex]{graphicx}
\DeclareGraphicsExtensions{.png}
\makeatletter
\renewcommand{\@biblabel}[1]{[#1]\hfill}
\makeatother

\usepackage[pdfencoding=auto]{hyperref}
\usepackage{bookmark}

\begin{document}
\setcounter{page}{1}

\begin{center}
{\LARGE \bf Metric spaces in chess and international chess pieces graph diameters\\[4mm]}
\vspace{3mm}

{\Large \bf Marco Rip\`a}
\vspace{3mm}

World Intelligence Network\\ 
Rome, Italy\\
e-mail: \url{marcokrt1984@yahoo.it}

\end{center}
\vspace{4mm}


\noindent {\bf Abstract:} This paper aims to study the graph radii and diameters induced by the $k$-dimensional versions of the well-known six international chess pieces on every finite $\{n \times n \times \dots \times n\} \subseteq \mathbb{Z}^k$ lattice since they originate as many interesting metric spaces for every proper pair $(n,k)$. For this purpose, we finally discuss a mathematically consistent generalization of all the planar FIDE chess pieces to an appropriate $k$-dimensional environment, finding the exact values of the graph radii and diameters of the $k$-rook, $k$-king, and $k$-bishop for every $k \geq 2$, as well as the corresponding values for the $3$-queen and $3$-pawn. We also give the diameter of the $3$-knight graph and provide tight bounds for the graph radii and diameters of the $k$-queen, $k$-knight, and $k$-pawn, holding for every $k \geq 4$.


\vspace{2mm}


\section{Introduction} \label{sec:Intr}
\sloppy Chess (i.e., ``international chess'' ruled by FIDE, the International Chess Federation \cite{1}) is the worldwide most popular board game with perfect information (together with Go \cite{2}). International chess is played by two opponents on an $8 \times 8$ checkerboard, where each player controls an army consisting of eight pawns ($P$), two knights ($N$), two bishops ($B$), two rooks ($R$), one queen ($Q$), and one king ($K$).
Each of the aforementioned $6$ kinds of chess pieces moves uniquely across the playing field.

On the other hand, in mathematics, a $\textit{metric space}$ is just a set of elements that has a notion of distance (i.e., a distance function called ``metric'') between every pair of elements belonging to the aforementioned set.

Thus, let $\mathbb{S}$ be a set of $n^k$ points ($n,k \in \mathbb{Z}^+ : \rm{P_1}, \rm{P_2},\dots,$ \hspace{-1mm}${\rm{P}}_{n^k-1},{\rm{P}}_{n^k} \in \mathbb{S}$) and let the integer-valued function $\delta$ be defined on pairs of elements of $\mathbb{S}$ so that $\mathbb{S}$ is called a metric space, having metric $\delta$, if all the following three conditions are satisfied \cite{12}.
\begin{enumerate}
\item $\delta(\rm{P_1}, \rm{P_2})\geq 0 \quad \forall \rm{P_1}, \rm{P_2} \in \mathbb{S} \quad \wedge \quad$ $\delta(\rm{P_1}, \rm{P_2})=0 \quad \textnormal{iff} \quad \rm{P_1}=\rm{P_2}$.
\item $\delta(\rm{P_1}, \rm{P_2})=$ $\delta(\rm{P_2}, \rm{P_1}) \quad \forall \rm{P_1}, \rm{P_2} \in \mathbb{S}$.
\item $\delta(\rm{P_1}, \rm{P_2}) \leq$ $\delta(\rm{P_1}, \rm{P_3})$ $+$ $\delta(\rm{P_3}, \rm{P_2}) \quad \forall \rm{P_1}, \rm{P_2}, \rm{P_3} \in \mathbb{S} \hspace{2mm}$ (the well-known triangle inequality).
\end{enumerate}

Now, as we usually do in standard geometry, we can easily extend the notion of the radius of a graph (by defining its diameter since the radius of a graph exists only if it has the diameter) to every metric space originated by any of the described chess pieces that verify the three fundamental conditions stated above \cite{10, 12}.

For the sake of clarity, from here on, let us consider any given, $k$-dimensional, hypercubic lattice defined (through the standard Cartesian product) as $C(n,k):=\{\{0,1,\dots,n-1\}\times \{0,1,\dots,n-1\} \times \cdots \times \{0,1,\dots,n-1\}\} \subset \mathbb{Z}^k$. We are interested in finding the diameter of the graphs (see Definition \ref{def1.2}) induced on $C(n,k)$ by the rook ($d_{n}^k(R)$), by the king ($d_{n}^k(K)$), by the queen ($d_{n}^k(Q)$), by the knight ($d_{n}^k(N)$), and we wish to do the same for the bishop and the pawn, bearing in mind that, in general, it is not allowed to let any edge go outside the AABB $[0,1,\dots,n-1] \times [0,1,\dots,n-1] \times \cdots \times [0,1,\dots,n-1]$ nor to consider any Steiner point whose coordinates include some non-integer values, of course.

All the chess pieces above (extended to multiple dimensions as in Millennium $3$D chess by William L. d'Agostino \cite[p.\ 227]{13}), with the sole exceptions of bishop and pawn, naturally induce a metric space on any $C(n,k) \subseteq \mathbb{Z}^k$ lattice such that $n \geq 4 \wedge k \geq 2$ (e.g., if $n = 3$ is given and we additionally assume that the knight can only make its well-known L-shaped move, then the knight distance between the vertex $(1, 1, \dots, 1)$ and any other vertex belonging to \linebreak $\{\{0,1,2\} \times \{0,1,2\} \times \cdots \times\{0,1,2\}\} \subset \mathbb{R}^k$ would be infinite since any knight move mandatorily changes by $2$ a coordinate and by $1$ another one so that, as we start from $(1, 1, \dots, 1)$, $1+2>2$ or $1-2<0$, and this would send the knight outside our virtual checkerboard at move $1$).

In this regard, let us point out that we are assuming the most obvious knight distance derived from the FIDE Handbook \cite{1}, but different metrics can also be considered to allow the knight to visit any vertex of given $k$-dimensional grids $C(3,k)$, including the central point \cite{22}. \linebreak As an example, we can set $k=5$ so that the Euclidean distance between $(0,0,0,0,0)$ and the center of $C(3,5)$ is equal to $\sqrt{1^2+2^2}$, which corresponds to the Euclidean distance between $(0,0)$ and $(2,1)$ in $\mathbb{R}^2$. In detail, we can anticipate the existence of open Euclidean knight's tours in $C(3,5)$, and consequently $C(3,5)$ is also a metric space for such a Euclidean knight. Furthermore, the same argument holds for every $C(2,k)$ such that $k \geq 7$ since it is possible to show that, if $k>6$ is given, the described Euclidean knight can also visit every vertex of $C(2,k)$ by performing a finite number of jumps (i.e., \cite[Theorem~4.1]{22} constructively proves that a closed Euclidean knight's tour is always possible in $C(2,k)$, as long as $k > 6$) of (Euclidean) length $\sqrt{5}$ units each, but this will not be the goal of the present paper because we are not going to circumvent any of the natural constraints arising from the FIDE Handbook.

\begin{figure}[H]
\begin{center}
\includegraphics[width=\linewidth]{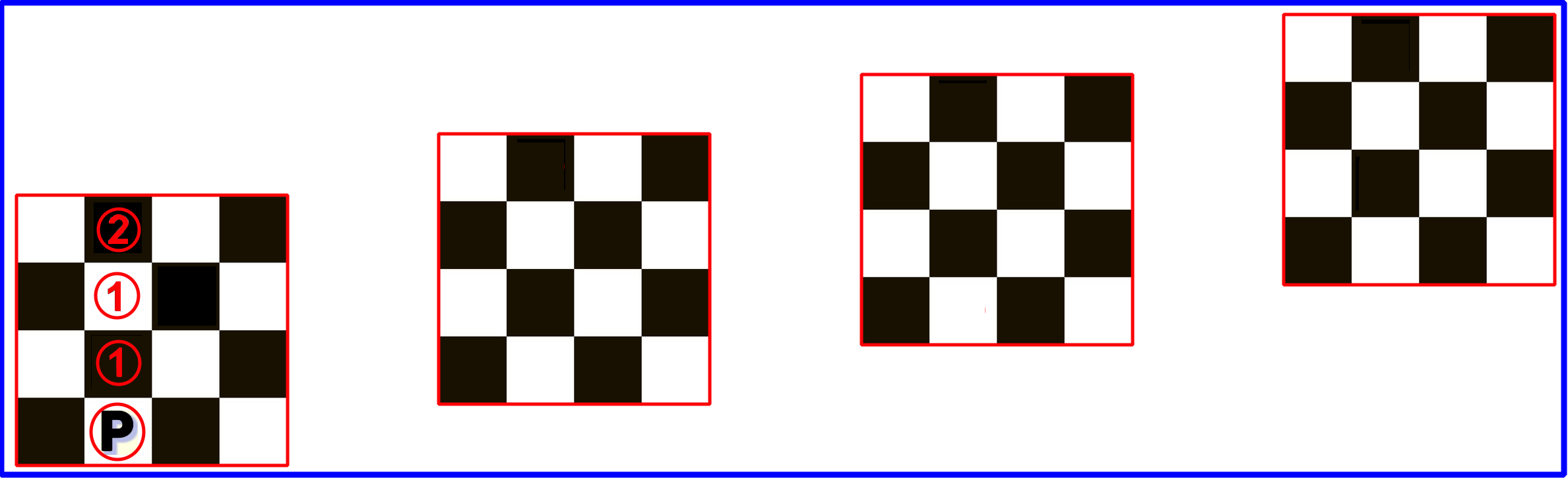}
\end{center}
\caption{A graphical representation of the sets $S_t(\overline{P};4,3;1,0,0)$ for $t \in \{0,1\}$ (e.g., if $t=2$, then $S_2(\overline{P};4,3;1,0,0)=\{(1,0,0), (1,1,0), (1,2,0), (1,3,0)\}$ follows).}
\label{fig:Figure_8}
\end{figure}

\begin{figure}[H]
\begin{center}
\includegraphics[width=\linewidth]{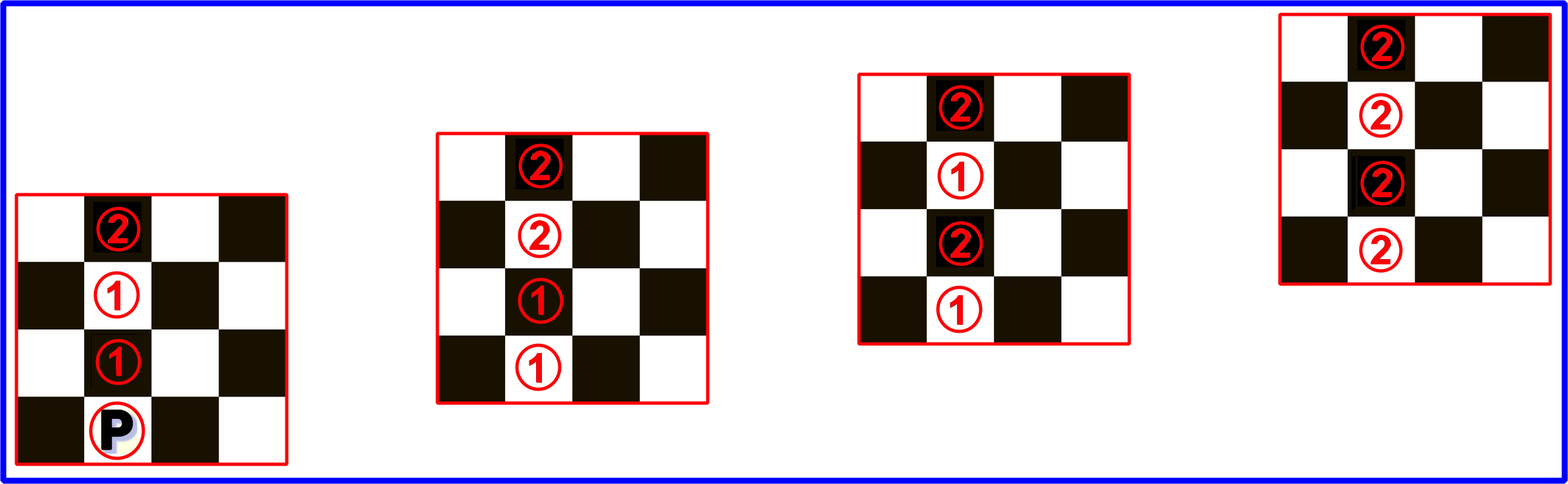}
\end{center}
\caption{A graphical representation of the sets $S_t(P;4,3;1,0,0)$ for any $t \in \{0,1\}$ (e.g., if $t=1$, then $S_1(P;4,3;1,0,0)=\{(1,0,0), (1,1,0), (1,2,0), (1,0,1), (1,1,1), (1,0,2), (1,2,2)\}$ follows), by using a $3$D checkerboard which is consistent with the Millennium $3$D chess environment.}
\label{fig:Figure_9}
\end{figure}

Thus, let $C(n,k)$ be such that $n \geq 4 \wedge k \geq 2$. It is easy to observe that the only two chess pieces that do not produce any metric space by their own are the bishop and the pawn (for nontrivial values of $n$, we will implicitly assume that every pawn promotes to a queen or knight after at most $n-2$ moves, see Figures \ref{fig:Figure_8}\&\ref{fig:Figure_9} at the end of the present section, and then we will provide a consistent metric for the $k$-pawn in Section \ref{sec:4}), but here we can simply state that the light-square bishop will always spend $2^k$ moves to reach any dark ``square'' of the ($k$-dimensional) board and vice versa.
Accordingly, let us introduce a few definitions by assuming that each chess piece induces a connected graph, covering every vertex of $C(n,k):=\underbrace{\{0,1, \dots, n-1\}\times\{0,1, \dots ,n-1\}\ \times \cdots \times\{0,1, \dots ,n-1\}}_\textrm{\textit{k}-times}$ for any pair $(n,k) : n \geq 4 \wedge {k \geq 2}$.

\begin{definition} \label{def1.1}
\sloppy By considering any vertex ${\rm{V}}:={\rm{V}}(n,k)$ belonging to the set $C(n,k)$, the eccentricity, $e_{n}^k(X;{\rm{V}})$, of the (connected) graph of any given chess piece $X$ is the distance between ${\rm{V}}$ and another vertex of $C(n,k)$ farthest from ${\rm{V}}$ itself.
\end{definition}

\begin{definition} \label{def1.2}
The diameter, $d_{n}^k(X)$, of the graph of any given chess piece $X$ is the maximum eccentricity among all the vertices belonging to the set $C(n,k)$ (see Definition \ref{def1.1}).
\end{definition}

In other words (by assuming $n>3$), Definition \ref{def1.2} suggests that the graph diameter of any chess piece is a positive integer that corresponds to the minimum number of moves needed by that piece, in the worst-case scenario, to go from one vertex, ${\rm{V_1}}$, of $C(n,k)$ to another vertex, ${\rm{V_2}}$, of the same set, where both of them form also the farthest pair of vertices of $\{\{0,1,\dots,n-1\}\times\{0,1,\dots,n-1\}\times \cdots \times\{0,1,\dots,n-1\}\}$ and vice versa (since a necessary but not sufficient condition for having a metric space is that $\delta({\rm{V_1}}, {\rm{V_2}})=\delta({\rm{V_2}}, {\rm{V_1}})$ holds for all ${\rm{V_1}}, {\rm{V_2}} \in C(n,k)$).

Thus, studying the graph diameters of every chess piece for as many pairs $(n,k)$ as we can is crucial.

\sloppy Unfortunately, from the FIDE Handbook, we can derive more than one reasonable way to extend the six chess pieces to the third dimension.
For this purpose, the remaining part of the present section provides a first, preliminary, definition of the move rules of the generalized versions of the FIDE chess pieces so that, in Section \ref{sec:2}, we will set the value of $k$ at $2$ for the chosen multidimensional versions of queen and knight, fully solving the problem for the extended king, rook, and bishop. On the other hand, the metric of the $k$-dimensional pawn, $\tilde{P}$, will be stated in Section \ref{sec:4} since it is not very easy to explain why $\delta_{n}^k(\tilde{P})({\rm{V_1}}, {\rm{V_2}}):= n-2+\max\{d_{n}^2(\tilde{Q}), d_{n}^2(\tilde{N}) \} +\min\{\delta_{n}^k(\tilde{Q})({\rm{V_1}}, {\rm{V_2}}),\delta_{n}^k(\tilde{N})({\rm{V_1}}, {\rm{V_2}}) \}$ would properly describe the pawn distance between the distinct vertices ${\rm{V_1}} \in C(n,k)$ and ${\rm{V_2}} \in C(n,k)$ for every $n > 3$.
For this reason, in the present section, we only introduce a pair of examples of pawn move rules that can be considered by $3$D chess players, even if this approach cannot produce any metric by itself.

In Section \ref{sec:3}, we will assume $k \geq 3$ and focus our analysis on the chosen generalizations of $Q$ and $N$. We then proceed by studying the radii and diameters of their graphs.

\sloppy Finally, we will adopt a strict mathematical approach, invoking the parity argument \cite{25} and applying retrograde check to the needed definitions of all the elements of the set $\{B, R, Q, K, N, P\}$, in order to describe the most consistent $k$-dimensional generalization of every planar chess piece introduced by the third article of the FIDE Handbook \cite{1}.

Hence, Section \ref{sec:4} will describe the above, being also devoted to providing nontrivial bounds for the radii and diameters of the corresponding graphs.

In the end, we strongly believe that, in a strict mathematical sense, the definitions given in Section \ref{sec:4} should be taken as the most straightforward and logical reference that extends the six FIDE chess pieces \cite{1} to multiple dimensions (we will refer to them as $k$-bishop or $\tilde{B}$, $k$-rook or $\tilde{R}$, $k$-queen or $\tilde{Q}$, and so forth), with the only partial exception of the Euclidean knight \cite{22}, which provides an alternative way to look at the classical $k$-knight \cite{17, 26, 27, 28} (see also the description provided in the third chapter of Watkins' book \cite{29})  when $k$ is above $4$.
\begin{definition} \label{def1.3}
The radius, $r_{n}^k(X)$, of the graph of any given chess piece $X$ is the minimum eccentricity among all the vertices belonging to the set $C(n,k)$ (see Definition \ref{def1.1}).
\end{definition}

At this point, it is worth mentioning a well-known graph theory theorem, stating that the diameter of any connected graph always lies between the radius, $r$, and its double \cite{11}. Since the knight, the rook, the queen, and the king (let us call $\hat{X}$ the selected chess piece belonging to $\{N,R,Q,K\}$) certainly produce connected graphs for every pair $(n,k) : n \geq 4 \wedge k \geq 2$, it follows that $r_{n}^k(\hat{X}) \leq d_{n}^k(\hat{X}) \leq 2 \cdot r_{n}^k(\hat{X})$.

\begin{definition} \label{def1.4}
\sloppy Let $X$ be a chess piece of the set $\{P, N, B, R, Q, K\}$. We indicate as $S_t(X;n,k;x_1,x_2,\dots,x_k) \subseteq C(n,k)$ the subset of all the vertices of $C(n,k)$ that can be reached in $t \in \{0,1,2,\dots,d_{n}^k(X)\}$ moves (or less) starting from the vertex identified by the $k$ Cartesian coordinates ($k$-tuple) $(x_1,x_2,\dots,x_k)$. In particular, if $t=1$, then $S_t(X;n,k;x_1,x_2,\dots,x_k)$ identifies the vertices of $C(n,k)$ that are zero or one $X$ move away from $(x_1,x_2,\dots,x_k)$, while $S_{d_{n}^k(X)+c}(X;n,k;x_1,x_2,\dots,x_k)=S_{d_{n}^k(X)}(X;n,k;x_1,x_2,\dots,x_k)=C(n,k)$ for all $c \in \mathbb{N}_0$.
\end{definition}

Now, for each of our six chess pieces, let us properly describe the mentioned subsets of all the vertices of $C(n,k)$ that are at most one move away from the starting vertex (i.e., $(x_1,x_2,\dots,x_k) \in C(n,k)$) so that we can better understand the problem of bounding $d_{n}^k(X)$ in higher dimensions (by Definition \ref{def1.4}, $X \in \{P, N, B, R, Q, K\}$ and $S_1(X;n,k;x_1,x_2 \dots,x_k) = S_0(X;n,k;x_1,x_2,\dots,x_k) \cup S_1(X;n,k;x_1,x_2,\dots,x_k)=\{(x_1,x_2,\dots,x_k)\} \cup S_1(X;n,k;x_1,x_2,\dots,x_k)$). We point out that the above is a necessary but not a sufficient condition for being sure that the two vertices are at a unit $X$ distance (e.g., they are at a unit knight distance from each other if $X=N$) in a metric space; this is clearly shown by the FIDE pawn on the regular $8 \times 8$ chessboard, as we can start the game by moving it from ${\rm{e}}2$ to ${\rm{e}}3$, not being allowed to put it back on ${\rm{e}}2$ at move $2$.

Let us start with the bishop case.

\begin{figure}[H]
\begin{center}
\includegraphics[width=\linewidth]{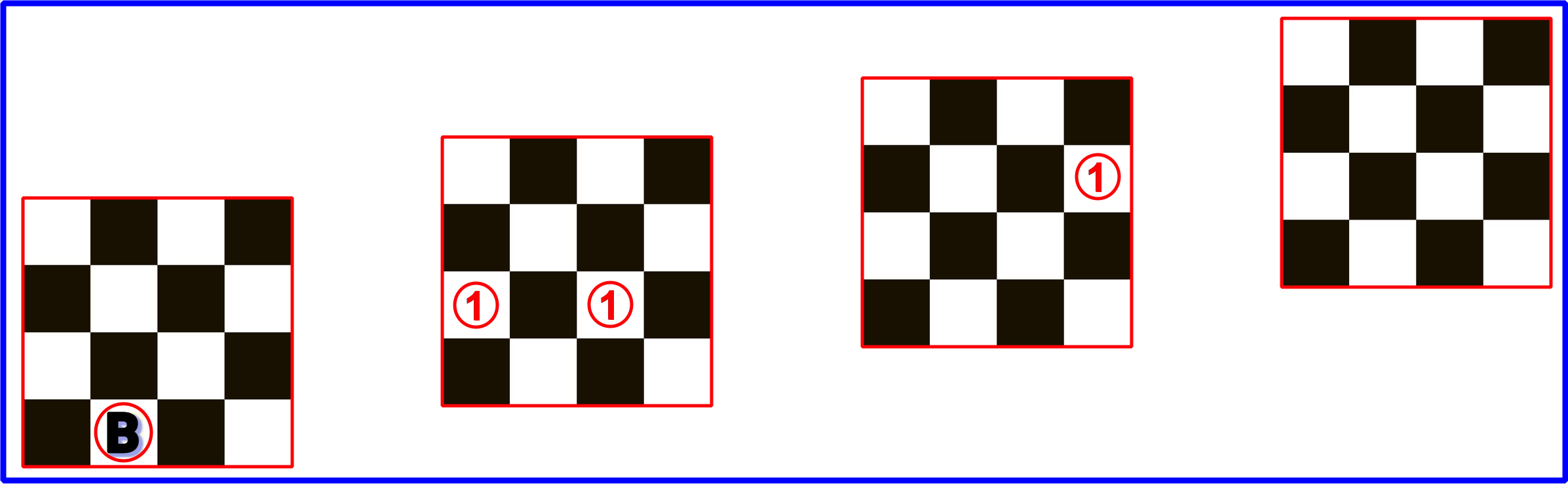}
\end{center}
\caption{A graphical representation of the set $S_1(\overline{B};4,3;1,0,0)=\{(1,0,0),(0,1,1), (2,1,1), (3,2,2)\}$.}
\label{fig:Figure_1}
\end{figure}

Since there is no unique way to extend the bishop move rule from $k=2$ to any $k \geq 3$, we can apply a sort of retrograde check to the other involved chess pieces to decide which move rule definition would lead to the most interesting configurations to be studied. We know that, in chess (see Reference \cite{1}, Articles 3.2 to 3.4), the queen is usually described as a piece that can alternatively move as a rook or as a bishop and so, if we define the bishop move rule as $S_1(\overline{B};n,k;x_1, x_2, \dots, x_k)=\{(x_1,x_2,\ldots,x_k)\} \cup \{(x_1+c_1, x_2+c_2, \dots, x_k+c_k) : \forall j \in \{1,2,\dots,k\}, \hspace{2mm} c_j \in \{-c,c\} \wedge (x_j+c_j) \in \{0,1,\dots,n-1\}\}$ (see Figure \ref{fig:Figure_1}), then the queen defined by $S_1(\overline{Q};n,k;x_1,x_2,\dots,x_k) := S_1(\overline{B};n,k;x_1,x_2,\dots,x_k) \cup S_1(\overline{R};n,k;x_1,x_2,\dots,x_k)$ could move from $(1,0,0)$ to $(2,1,1)$ in one move on $C(4,3)$, while reaching $(2,1,0)$ from the same starting vertex would require two moves.

\begin{figure}[H]
\begin{center}
\includegraphics[width=\linewidth]{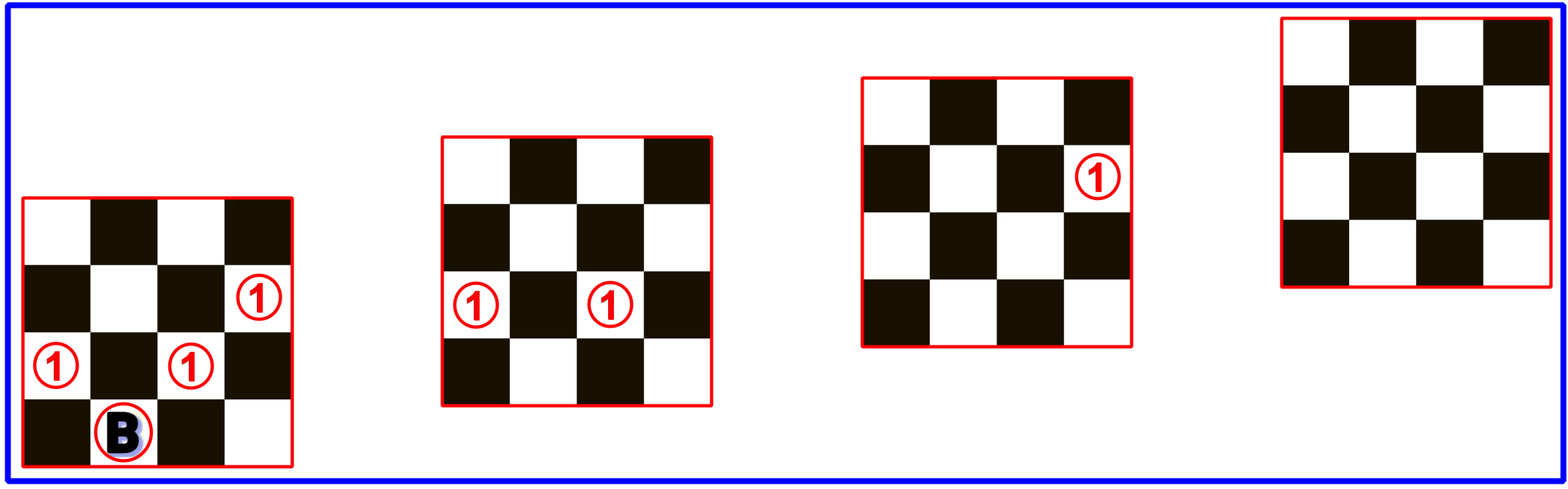}
\end{center}
\caption{A graphical representation of the set ${S_1(B;4,3;1,0,0)}=\{(1,0,0), (0,1,0), (2,1,0), (3,2,0), (0,1,1), (2,1,1), (3,2,2)\}$, \linebreak by using a $3$D checkerboard which is consistent with the Millennium $3$D chess environment.}
\label{fig:Figure_2}
\end{figure}

Therefore, it would seem reasonable to state the (generalized) bishop move rule as follows (see Figure \ref{fig:Figure_2}).

Let $x_1, x_2, \dots, x_k \in \{0,1,\dots,n-1\}$ be given. Then, $S_0(B;n,k;x_1,x_2,\dots,x_k):=\{(x_1, x_2, \dots, x_k)\}$ and $S_1(B;n,k;x_1, x_2, \dots, x_k)= S_0(B;n,k;x_1,x_2,\dots,x_k) \cup \{(x_1+c_1, x_2+c_2, \dots, x_k+c_k) : \forall j \in \{1,2,\dots,k\}, \hspace{2mm} (c_j \in \{-c,0,c\} \hspace{3mm} \wedge c_j \neq 0 \quad \textnormal{if} \quad j < 3) \wedge (x_j+c_j) \in \{0,1,\dots,n-1\}\}$ (e.g., if we take as our starting point the vertex $(1,0,0) \in C(4,3)$, it follows that $S_1(B;4,3;1,0,0) = \{(1,0,0), (0,1,0), (2,1,0), (3,2,0), (0,1,1), (2,1,1), (3,2,2)\}$).

\begin{figure}[H]
\begin{center}
\includegraphics[width=\linewidth]{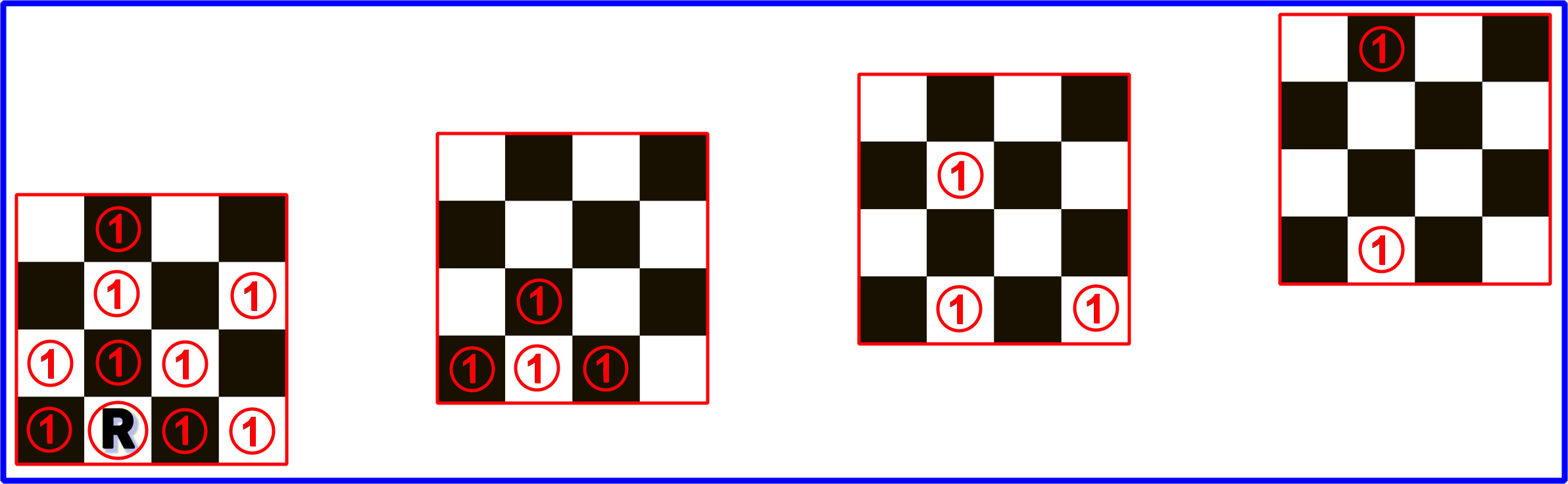}
\end{center}
\caption{A graphical representation of the set $S_1(R^{\star};4,3;1,0,0)=\{(1,0,0), (0,0,0), (2,0,0), (3,0,0), (0,1,0), (1,1,0),  (2,1,0), (1,2,0), (3,2,0), (1,3,0), \linebreak (0,0,1), (1,0,1), (2,0,1),$ $(1,1,1), (1,0,2), (3,0,2), (1,2,2), (1,0,3), (1,3,3)\}$.}
\label{fig:Figure_3bis}
\end{figure}

\begin{figure}[H]
\begin{center}
\includegraphics[width=\linewidth]{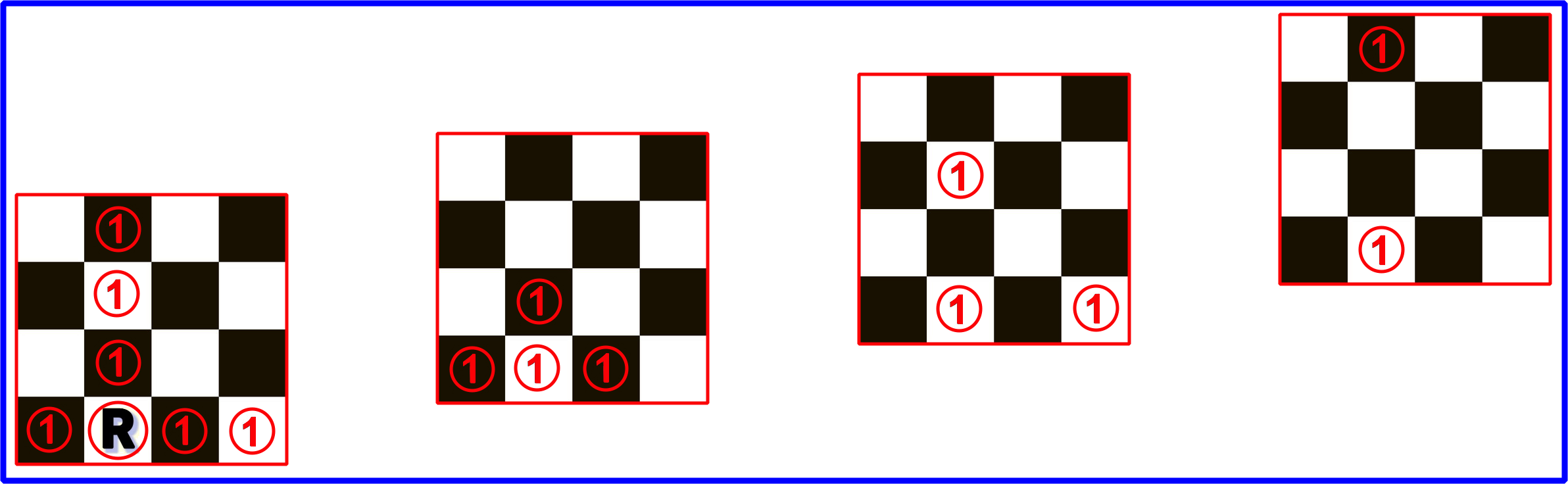}
\end{center}
\caption{A graphical representation of the set $S_1(R;4,3;1,0,0)=\{(1,0,0), (0,0,0), (2,0,0), (3,0,0), (1,1,0), (1,2,0), (1,3,0), (0,0,1), (1,0,1), (2,0,1),$ $(1,1,1), (1,0,2), (3,0,2), (1,2,2), (1,0,3), (1,3,3)\}$, by using a $3$D checkerboard which is consistent with the Millennium $3$D chess environment.}
\label{fig:Figure_3}
\end{figure}

\begin{figure}[H]
\begin{center}
\includegraphics[width=\linewidth]{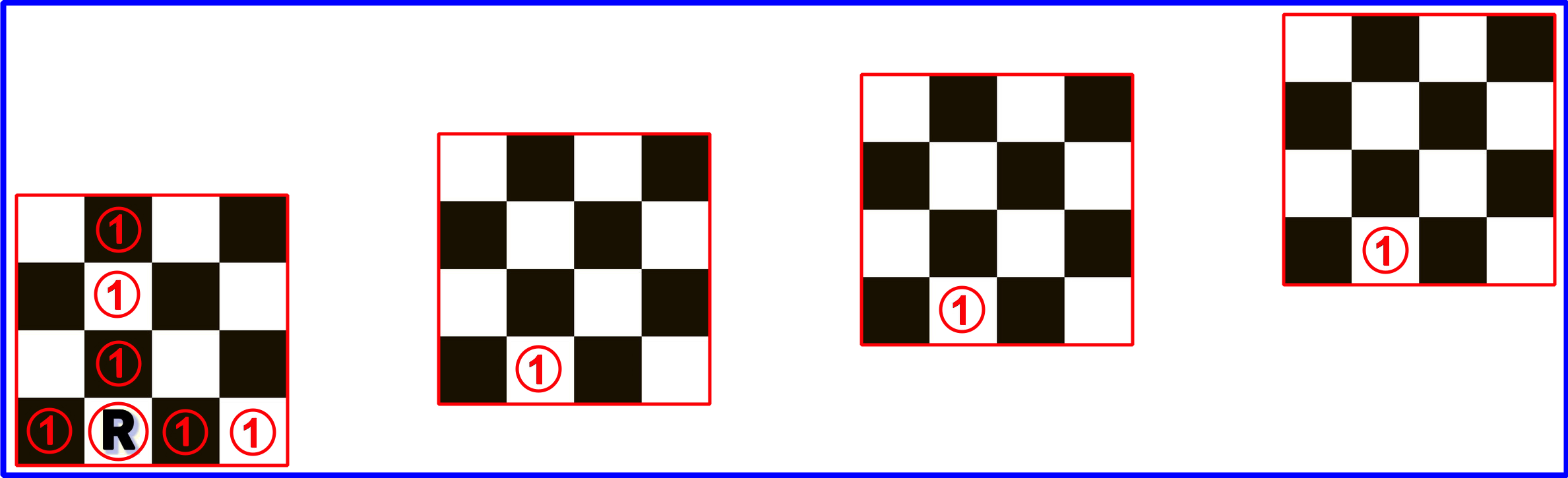}
\end{center}
\caption{A graphical representation of the set $S_1(\overline{R};4,3;1,0,0)=\{(1,0,0), (0,0,0), (2,0,0), (3,0,0), (1,1,0), (1,2,0), (1,3,0), (1,0,1), (1,0,2), (1,0,3)\}$.}
\label{fig:Figure_4}
\end{figure}

Since $S_t(X;n,k;x_1,x_2,\dots,x_k) = S_{t-1}(X;n,k;x_1,x_2,\dots,x_k) \cup S_t(X;n,k;x_1,x_2,\dots,x_k)$ trivially holds for any $t \in \mathbb{Z}^+$, by using a very similar argument to that previously considered for choosing $S_1(B;n,k;x_1, x_2, \dots, x_k)$ over $S_1(\overline{B};n,k;x_1, x_2, \dots, x_k)$, we could define the extended rook, $R^{\star}$, as shown in Figure \ref{fig:Figure_3bis} (i.e., a $k$-dimensional rook $R^{\star}$ can directly move from the vertex ${\rm{V_1}} \equiv (x_1, x_2, \dots, x_k) \in C(n,k)$ to another vertex, \linebreak${\rm{V_2}} \equiv (y_1, y_2,\dots,y_k)$, if there exists some constant $c \in \{0,1,\ldots,n-1\}$ such that \linebreak$|x_j - y_j| \in \{0,c\}$ for each $j=1, 2, \dots, k$; additionally, such a rook cannot change all its $k$ Cartesian coordinates at a time like the bishops shown in Figures \ref{fig:Figure_1}\&\ref{fig:Figure_2} so that there must be at least one integer $j$ such that $x_j = y_j$).
Accordingly, let $n, k \in \mathbb{N}-\{0,1\}$ and assume $S_0(X;n,k;x_1, x_2, \dots, x_k)=\{(x_1, x_2, \dots, x_k)\}$. If $X:=R^{\star}$ represents such a powerful rook, we can write that $S_1(R^{\star};n,k;x_1,x_2,\dots,x_k)=\{(x_1,x_2,\ldots,x_k)\} \cup \linebreak \{(y_1,y_2,\ldots,y_k)\in C(n,k) : \exists c\in\{1,2,\ldots,n-1\},\ \forall j\in\{1,2,\ldots,k\},\ |y_j-x_j|\in\{0,c\},\quad 1\leq |\{j:y_j\neq x_j\}|\leq k-1\}$. Unfortunately, the mentioned move rule would not describe the Millennium $3$D chess rook, $R$, so we are forced to state that $S_1(R;n,k;x_1,x_2,\dots,x_k) :=\{(x_1,x_2,\ldots,x_k)\} \cup \bigl(S_1(R^{\star};n,k;x_1,x_2,\dots,x_k) - S_1(B;n,k;x_1,x_2,\dots,x_k)\bigr)$ (or, at least, we cannot deny the subset relation $S_1(R;n,k;x_1,x_2,\dots,x_k) \subseteq \{(x_1,x_2,\ldots,x_k)\} \cup \bigl(S_1(R^{\star};n,k;x_1,x_2,\dots,x_k) - S_1(B;n,k;x_1,x_2,\dots,x_k)\bigr)$, see Figure \ref{fig:Figure_3}). In this section, we skip the most straightforward move rule for an extended rook, given by the subset $S_1(\overline{R};n,k;x_1,x_2,\dots,x_k) = \bigcup_{1 \leq j \leq k} \{ (x_1, \ldots, x_{j-1}, c, x_{j+1}, \ldots, x_k) \hspace{2mm} : \hspace{2mm} c \in \{0, 1, \ldots, n-1\} \}$, as shown in Figure \ref{fig:Figure_4}, since we will properly discuss it in Section \ref{sec:4}, when we will select the most reliable $k$-dimensional generalization for each of the six international chess pieces.

\begin{figure}[H]
\begin{center}
\includegraphics[width=\linewidth]{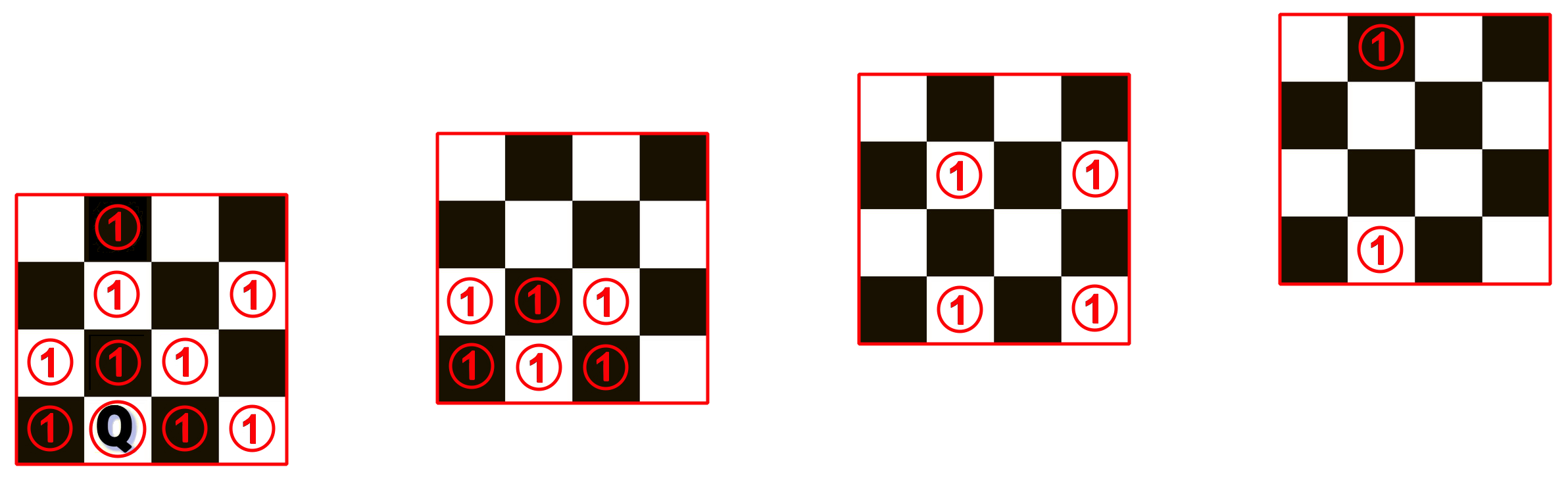}
\end{center}
\caption{A graphical representation of the set $S_1(Q;4,3;1,0,0) = \{(1,0,0), (0,0,0),$ $(2,0,0), (3,0,0), (1,1,0), (1,2,0), (1,3,0), (0,0,1), (1,0,1), (2,0,1),$ $(1,1,1), (1,0,2),$ $(3,0,2), (1,2,2), (1,0,3), (0,1,0), (2,1,0), (3,2,0), (0,1,1), (2,1,1), (3,2,2), (1,3,3)\}$,\linebreak by using a $3$D checkerboard which is consistent with the Millennium $3$D chess environment.}
\label{fig:Figure_5}
\end{figure}

The queen case follows by taking the union of the moves of $R$ and $B$ (equivalently, of $R^{\star}$ and $B$).
Let $n, k \in \mathbb{N}-\{0,1\}$. If $X:=Q$, then $S_1(Q;n,k;x_1,x_2,\dots,x_k) = S_1(B;n,k;x_1,x_2,\dots,x_k) \cup S_1(R;n,k;x_1,x_2,\dots,x_k)$ (see Figure \ref{fig:Figure_5}).

\begin{figure}[H]
\begin{center}
\includegraphics[width=\linewidth]{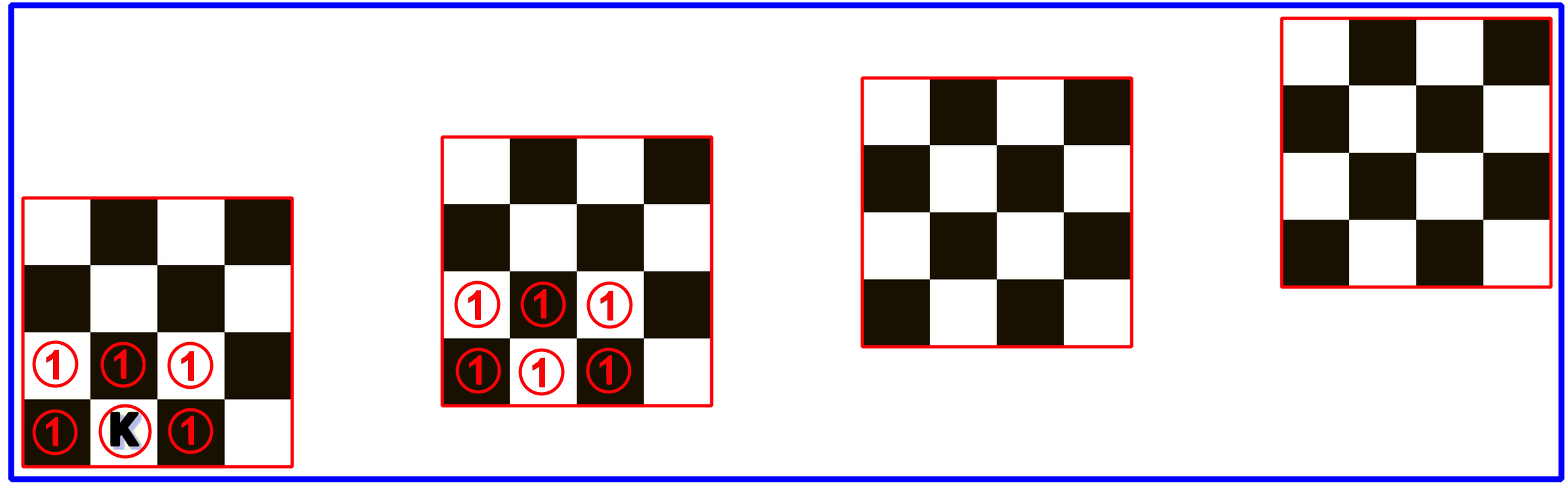}
\end{center}
\caption{A graphical representation of the set $S_1(K;4,3;1,0,0)=\{(1,0,0), (0,0,0), (2,0,0), (0,1,0), (1,1,0), (2,1,0), (0,0,1), (1,0,1), (2,0,1), (0,1,1),$ $(1,1,1), (2,1,1)\}$, by using a $3$D checkerboard which is consistent with the Millennium $3$D chess environment.}
\label{fig:Figure_6}
\end{figure}

The king case is described by \v{C}heby\v{s}\"ev's distance $\delta_{Chebyshev}((x_1,\ldots,x_k),(y_1,\ldots,y_k)):=\max_j\{|x_j-y_j|, \hspace{2mm} j=1,\ldots,k\}$ \cite{50}. Accordingly, as usual, let $n, k \in \mathbb{N}-\{0,1\}$. If $X:=K$, then $S_1(K;n,k;x_1,x_2,\dots,x_k) = \{(x_1+c_1,x_2+c_2,\dots,x_k+c_k) : c_1, c_2, \dots, c_k \in \{-1,0,1\} \hspace{2mm} \wedge \hspace{2mm} (x_j+c_j) \in \{0,1,\dots,n-1\} \hspace{2mm} \forall j \in \{1,2,\dots,k \}\}$ (see Figure \ref{fig:Figure_6}).

\begin{figure}[H]
\begin{center}
\includegraphics[width=\linewidth]{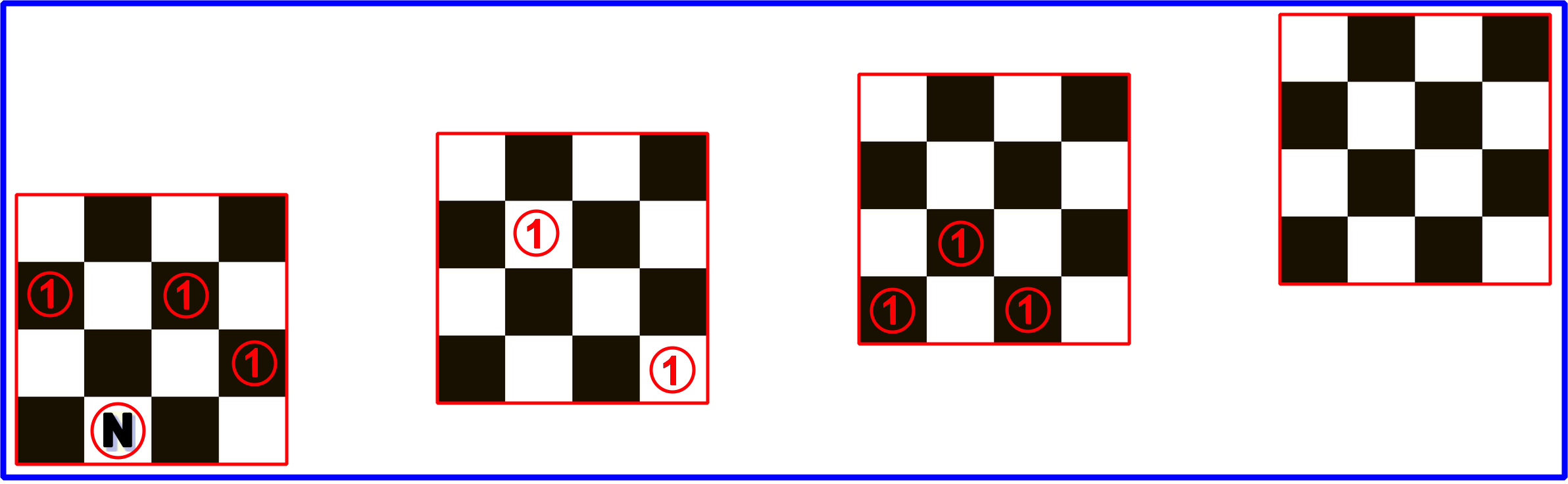}
\end{center}
\caption{A graphical representation of the set $S_1(N;4,3;1,0,0)=\{(1,0,0), (0,2,0), (2,2,0), (3,1,0), (1,2,1),$ $(3,0,1), (0,0,2), (1,1,2), (2,0,2)\}$, by using a $3$D checkerboard which is consistent with the Millennium $3$D chess environment.}
\label{fig:Figure_7}
\end{figure}

The knight case, which is also interesting for the related open knight's tour problem and closed knight's tour problem \cite{30} in higher dimensions \cite{22} (the latter solved for rectangular $k$-dimensional boards by J.\ Erde, B.\ Gol{\'e}nia, and S.\ Gol{\'e}nia in 2012 \cite[Thm.\ 3, pp.\ 3--4]{17}).
\linebreak
Let $n, k \in \mathbb{N}-\{0,1\}$. If $X:=N$, then $S_1(N;n,k;x_1, x_2, \dots, x_k)=\{(x_1, x_2, \dots, x_k)\} \cup \linebreak \{(x_1+c_1, x_2+c_2, \dots, x_k+c_k) \hspace{-0.5mm} : \hspace{-0.5mm} \forall j \hspace{-0.5mm} \in \hspace{-0.5mm} \{1,2,\dots,k\} \bigl( \sum_{j=1}^{k} |c_j| = 3 \wedge |c_j| \hspace{-0.5mm} \in \hspace{-0.5mm} \{0,1,2\} \hspace{-0.5mm} \wedge \hspace{-0.5mm} \exists! j \hspace{-0.5mm} : \hspace{-0.5mm}
|c_j|=1 \bigr) \wedge \linebreak (x_j+c_j) \in \{0,1,\dots,n-1 \}\}$ (see Figure \ref{fig:Figure_7}).

\sloppy About the pawn case (for any given pair $(n,k) : n,k \geq 2$), since there are no other pieces placed on our $k$-dimensional board, for our purposes, we could simply choose to consider the subset $S_1(\overline{P};n,k;x_1,x_2,\dots,x_k) =\{(x_1, x_2, \dots, x_k)\} \cup \{(x_1, x_2+c, x_3, \dots, x_k) : c \in \{1,2\} \wedge \left(x_1, x_2+c, x_3, \dots, x_k \in \{0,1,\dots,n-1 \}\right)\}$ (see Figure \ref{fig:Figure_8}) of ${S_1(P;n,k;x_1,x_2,\dots,x_k)} = \{(x_1, x_2+c_2, x_3+c_3, \dots, x_k+c_k) : ((c_j \in \{0,1\} \hspace{2mm} \forall j \in \{2,3, \dots, k\}) \vee (c_j \in \{0,2\} \hspace{2mm} \forall j \in \{2,3, \dots, k\})) \wedge (x_j+c_j) \in \{0,1, \dots, n-1\}\}$ (see Figure \ref{fig:Figure_9}). On the other hand, after the first move, it is well known that any pawn can only move forward by one chessboard unit at a time (see Article 3.7 of the FIDE Handbook).

Thus, if $n \geq 3$, we may only say that $S_1(\overline{P};n,k;x_1,x_2,\dots,x_k) =\{(x_1, x_2, \dots, x_k)\} \cup \{(x_1, x_2+c, x_3, \dots, x_k) : c \in \{1,2\} \wedge x_1, x_2+c, x_3, \dots, x_k \in \{0,1,\dots,n-1 \}\}$ and be satisfied (Figure \ref{fig:Figure_8}).

Anyway, when a pawn finally reaches a vertex in position $(x_1, n-1, x_3, \dots, x_k)$ it automatically promotes to a queen or knight, and it does not really matter if $S_t(\overline{P};n,k;x_1,x_2,x_3,\dots,x_k)$ is a too restrictive definition for some popular chess variants as d'Agostino's Millennium $3$D chess (see Figure \ref{fig:Figure_9}) \cite{13}, since it cannot make any sense to write $d_{n}^k(P)$ until we can properly define the $k$-pawn metric, with its symmetry axiom such that $\delta_{n}^k(P)({\rm{V_1}}, {\rm{V_2}})=\delta_{n}^k(P)({\rm{V_2}}, {\rm{V_1}})$ holds for any ${\rm{V_1}}$ and ${\rm{V_2}}$ belonging to $C(n,k)$, and this will be accomplished in Section \ref{sec:4}.

Lastly, it is trivial to point out that $S_0(K)=S_0(R)=S_0(B)=S_0(Q)=S_0(N)=S_0(P)=S_1(B) \cap S_1(R)=S_1(\overline{B}) \cap S_1(\overline{R})=\{(x_1, x_2, \dots, x_k)\}$, and this completes the introduction to our six $k$-dimensional chess pieces.


\section{Graph diameters on planar checkerboards} \label{sec:2}

It is quite reasonable to start our analysis from the easiest scenario, so let us consider planar graphs first (knowing that if $n=1$, then the diameter will always be equal to $0$, for any $k \geq 2$, regardless of which chess piece we are willing to consider).

Then, it is easy to note that the maximum number of moves for a FIDE king corresponds to half of the Manhattan distance between the point $(0,0)$ and $(n-1,n-1)$ so that $d_{n}^2(K)=n-1$ for any $n$, given the fact that, by definition, the $2$-king can change by one both the coordinates at the same time (e.g., it can go from $(1,1)$ to $(2,2)$ in only one move) or just increase/decrease by one a single coordinate (it is worth to note that, as we have anticipated in Section \ref{sec:Intr}, the $2$-king metric corresponds to the well-known Chebyshev distance, $L_{\infty}$ metric, or chessboard metric \cite{50}).
The same argument holds for every $k$ above $2$ since the $k$-king, in our construction (see Section \ref{sec:Intr}), can increase/decrease by one the value of $1$ to $k$ out of $k$ coordinates at the same time (i.e., by spending only one move). Thus, we will avoid repeating ourselves in Section \ref{sec:4} and already state that $d_{n}^k(K)=n-1$ holds for any $k$.

The rook (see Figures \ref{fig:Figure_3}\&\ref{fig:Figure_4}), by construction, can always change the value of one coordinate (by $1$ to $n-1$, depending on the starting vertex) at a time so that, in general, $\nexists n,k \in \mathbb{Z}^+ : d_{n}^k(\overline{R}) > k$. Hence, $d_{n}^k(\overline{R}) \leq k$ holds for any pair of positive integers $(n,k)$.
Then, let us consider that, $\forall t,n-1,k-2 \in \mathbb{N}_0$, $S_t(\overline{R};n,k;x_1, x_2, \dots, x_k) \subseteq {S_t(R;n,k;x_1, x_2, \dots, x_k)}$. It follows that $d_{n}^k(\overline{R}) \geq d_{n}^k({R})$. Thus, we have proven that $d_{n}^k(R) \leq k$.
Now, let us observe that, for any $n,k \in \mathbb{N}-\{0,1\}$, the $\overline{R}$-distance between the vertex $(0,0,\dots ,0) \in C(n,k)$ and $(1,1,\dots ,1) \in C(n,k)$ cannot be less than $k$ (i.e., $\delta_{n}^k(\overline{R})((0,0,\ldots,0),(1,1,\ldots,1))=k$, for every pair $(n,k)$ as above) so that $d_{n}^k(\overline{R}) \geq e_{n}^k({\overline{R};0,0,\dots,0}) \geq k$ implies $d_{n}^k(\overline{R}) = r_{n}^k(\overline{R}) = k$ ($r_{n}^k(\overline{R}) = k$ since $\delta_{n}^k(\overline{R})((x_1,x_2,\ldots,x_k),(x_1+1,x_2+1,\ldots,x_k+1))=k$ as long as $x_1,x_2,\ldots,x_k \in \linebreak \{0,1,\ldots,n-2\}$).

On the other hand, we can verify that for each value of $k \geq 3$, there exists a positive integer, $m$, that satisfies $\bigl(r_{m}^k(R)<d_{m}^k(R) \bigr) \wedge \bigl(r_{m+c}^k(R)=d_{m+c}^k(R) \bigr)$ for all $c \in \mathbb{Z}^+$ (e.g., if $k=3$, then $m=3$, since $r_{3}^3(R)=2<d_{3}^3(R)=3$ and $r_{n}^3(R)=d_{n}^3(R)=3$ for every $n \geq 4$), while $e_{3}^k(R;1,1,\ldots,1)=2$ implies that $e_{3}^k(R;1,1,\ldots,1)<e_{3}^k(R;0,0,\ldots,0)$ since $\delta_{3}^3(R)((0,0,0),(2,2,1))=3$ (so that $e_{3}^k(R;0,0,\ldots,0)$ cannot be smaller than $3$ as long as $k \geq 3$ is given).

Therefore, there is no reason to study the rook in Section \ref{sec:3}.

By definition, the queen alternatively moves as a rook or as a bishop, so the diameter of its graph cannot be greater than $\min\{d_{n}^k(R), d_{n}^k(B)\}$. Since $\min\{d_{n}^k(R), d_{n}^k(B)\}=d_{n}^k(R)$, it follows that $1 \leq d_{n}^k(Q) \leq k$ for all $n,k \in \mathbb{N}-\{0,1\}$. Furthermore, it is not difficult to see that, for any $k \geq 2$, $d_{2}^k(Q)=1$ and we can also check by brute force the third line of Equation (\ref{eq53}) so that its last line follows
\begin{equation}\label{eq53}
d_{n}^3(Q)=
 \begin{cases}
 0, & \text{iff $ n=1$;} \\
 1, & \text{iff $ n=2$;} \\
 2, & \text{iff $ n \in \{3,4,5\}$;} \\
 3, & \text{iff $ n \in \mathbb{N}-\{0,1,2,3,4,5\}$.}
 \end{cases}
\end{equation}
(see Figure \ref{fig:Figure_11}).

Trivially, we observe that $d_{1}^2(Q)=0$, $d_{2}^2(Q)=1$, and $d_{n \geq 3}^2(Q)=2$ since a queen cannot go from the vertex $(0,0)$ to $(1,2)$ in only one move, nor go from $(0,0,0)$ to $(5,4,2)$ by spending less than $3$ moves, and being $d_{n}^2(R)=2$ true for any $n \neq 1$, it follows that $d_{n}^2(Q)<d_{n}^2(R)$ if and only if $n=2$. Similarly, as we assume $k=3$, from $n > 2 \Rightarrow d_{n}^3(R)=3$ (which easily follows by observing that $\delta_{3}^3(R)((0,0,0), {(2,2,1)})>2$) and the consideration above, we deduce that $6$ is the smallest integer $n^{\star}\in \mathbb{Z}^+$ such that $d_{n}^3(Q)=d_{n}^3(R)=d_{n}^3(\overline{R})=3$ is true for every $n \geq n^{\star}$, but this pattern let us also conjecture that $\forall k \in \mathbb{N}-\{0,1\}$ $\exists \hspace{1mm} \mathbb{Z}^+ \hspace{-1mm} \ni n^{\star}:=n^{\star}(k) : d_{n^{\star}+q}^k(Q)=d_{n^{\star}}^k(Q)=k \hspace{1mm}, \forall q \in \mathbb{Z}^+$. In detail, $ n^{\star}=3$ clearly holds for $k=2$, while $n^{\star}=6$ is the solution for $k=3$, but what about the values of $n^{\star}$ (if any) satisfying $d_{n^{\star}}^k(Q)=k$ for $k=4,5,6,\ldots$?

\begin{figure}[H]
\begin{center}
\includegraphics[width=\linewidth]{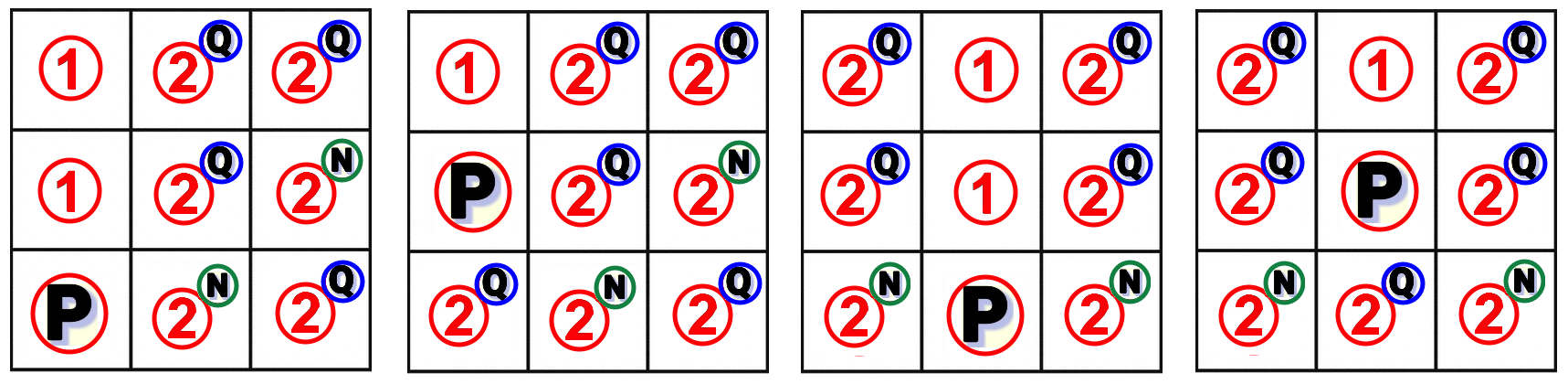}
\end{center}
\caption{Graphical proof that a pawn can reach any square on the $3 \times 3$ board in $2$ moves (omitting the trivial configurations where the starting point is the vertex $(0,2)$ or $(1,2)$).}
\label{fig:Figure_15}
\end{figure}

Here, according to the standard FIDE rules \cite{1}, we assume that a pawn can only move forward (so it cannot induce any metric on $C(n,k)$) and, if it stands in the first rank, it can alternatively increase that coordinate by one or by two, while all the other coordinates remain unchanged until it promotes (to a queen or a knight). Thus, a pawn spends no more than $n-2$ moves to promote and, since it is free to promote to a queen, in general, it can reach any given vertex of $C(n,k)$ in (at most) $d_{n}^k(Q)+n-2$ moves. On the other hand, it only takes $2$ moves (or less) for a pawn to move between any two vertices of $C(3,2)$ since it can go from any vertex on the third rank to any other vertex in only one move by choosing to promote to a knight or a queen, depending on the final square to be reached (see Figure~\ref{fig:Figure_15}).

When we take into account the knight (the most infamous piece on the board) and its graph diameter, we can find the value of $d_{n}^2(N)$ by just taking a look at the OEIS sequence A232007 \cite{15}, while the general $k=3$ case is fully covered by the OEIS sequence A359740 \cite{16} since the formal proof that $n \geq 4 \Rightarrow d_{n}^3(N) = n$ is by Pontus von Br\"omssen and Michael S. Branicky and it has been posted as a comment of the aforementioned OEIS sequence.
Thus, the only remaining cases are those such that $k \geq 4$, and this will be the topic of Subsection \ref{sec:sub4.2}.

Finally, the sixth chess piece to be analyzed is the bishop and, in order to get a finite value of $d_{n}^2(B)$ for every $n \in \mathbb{Z}^+$, we simply verify the triangle inequality by imposing that our light-square bishop can reach any dark square by spending a total of $2^k$ moves, and, symmetrically, we state that the dark-square bishop can reach any given light square after exactly $2^k$ moves.

In order to realize how to color every $k$-dimensional checkerboard, we need to understand that this directly follows from the definition of $S_1(B;n,k;x_1, x_2, \dots, x_k)$ since any light-square bishop cannot reach a dark square faster than any light square (and vice versa).

\begin{figure}[H]
\begin{center}
\includegraphics[width=\linewidth]{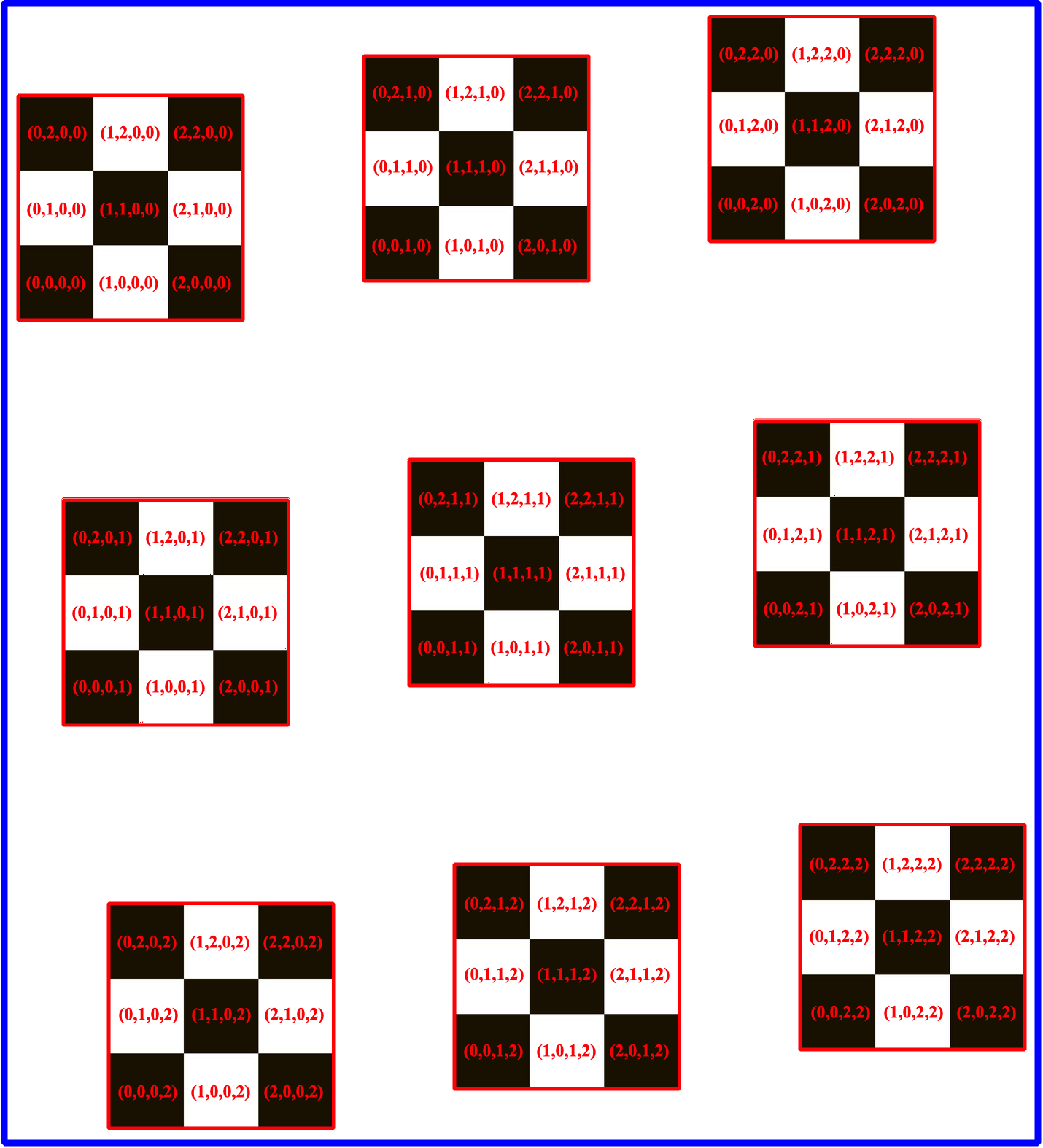}
\end{center}
\caption{Coloring $C(3,4)$ according to the stated, provisional, multidimensional bishop move~rule.}
\label{fig:Figure_10}
\end{figure}

As an example, the aforementioned constraint does not allow us to state that the coloring rule ``${\rm{V}} \equiv (x_1, x_2, \dots, x_k)$ belongs to the dark squares set if and only if $\sum_{j=1}^{k} x_j$ is an even integer, otherwise ${\rm{V}}$ is an element of the light squares subset of $C(n,k)$'' is consistent with our previous definition of bishop (and queen). So, Figure \ref{fig:Figure_10} simply describes one acceptable environment for the next section (while, in Section \ref{sec:4}, we will finally provide a better solution), representing the $n=3 \wedge k=4$ checkerboard forced by our provisional definition of bishop in three and more dimensions.

\section{Diameters in three dimensions and above} \label{sec:3}

In this section we discuss the values of $d_{n}^k$ for the remaining chess pieces, those we were not able to determine in Section \ref{sec:2}. They are the queen, the knight, and the pawn (since $d_{n}^k(P)$ cannot be defined for the reason discussed in the previous sections, we will consider only the minimum number of moves necessary to go from the starting vertex of $C(n,k)$ to the farthest one of the same set, by assuming the worst-case scenario). 

For this purpose, we observe that $d_{n}^3(P)$ will be (tightly) bounded by knowing $d_{n}^3(Q)$ (since $d_{n}^3(N)$ is fully described by the OEIS sequence A359740 and it is much greater than $d_{n}^3(Q)$ if $n \neq 1$), so our main goal is to find $d_{n}^3(Q)$ for any nontrivial value of $n$ (i.e., $n \geq 4$) and then we are finally free to follow the knight's path in higher dimensions.

\vspace{3mm}

Since we have previously shown that
\begin{equation}\label{eq54}
d_{n}^2(Q)=
 \begin{cases}
 0, & \text{iff $n=1$;} \\
 1, & \text{iff $n=2$;} \\
 2, & \text{iff $n \in \mathbb{N}-\{0,1,2\}$.} \\
 \end{cases}
\end{equation}
and also that
\begin{equation}\label{eq55}
d_{n}^3(Q)=
 \begin{cases}
 0, & \text{iff $n=1$;} \\
 1, & \text{iff $n=2$;} \\
 2, & \text{iff $n \in \{3,4,5\}$;} \\
 3, & \text{iff $n \in \mathbb{N}-\{0,1,2,3,4,5\}$.} \\
 \end{cases} \textnormal{\quad,}
\end{equation}
we can be interested in exploring the smallest values of $d_{n}^4(Q)$ to guess the next value of the integer sequence $a(k):=\min\{n \in \mathbb{Z}^+ :d_{n}^k(Q)=k\}$, which is $a(k)=3, 6, ?$ ($k \in \mathbb{N}-\{0,1\}$).

\begin{figure}[H]
\begin{center}
\includegraphics[width=\linewidth]{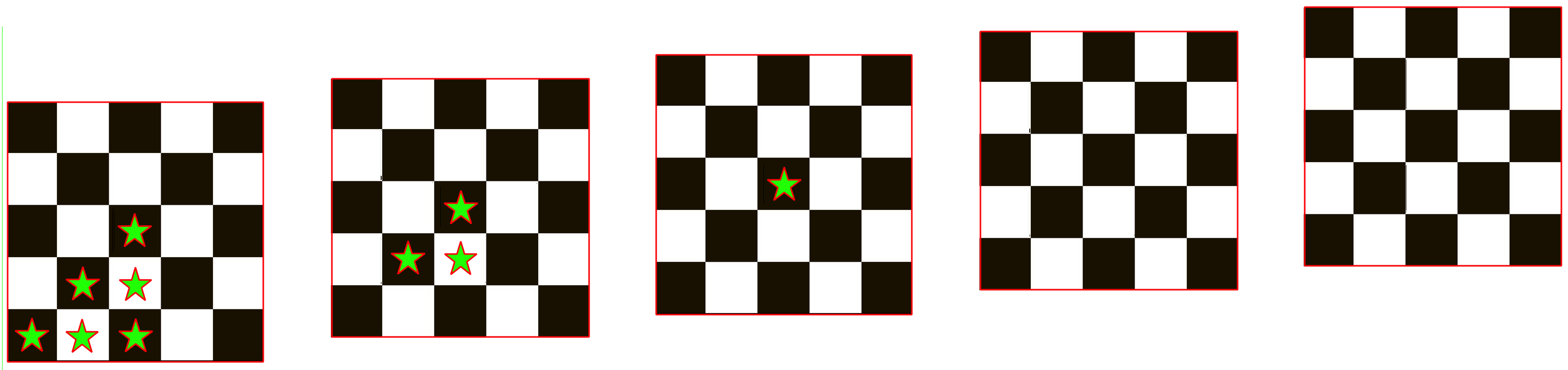}
\end{center}
\caption{Let $n=5$ (or even $n=6$), $k=3$, and assume that the chess piece $X$ is selected \linebreak (i.e., $X:=Q$). Then, by calculating the eccentricity of all the $10$ vertices marked with a green star, we automatically solve by brute force the problem of finding $d_{5}^3(X)=\max\{e_{5}^3(X;0,0,0), e_{5}^3(X;1,0,0), e_{5}^3(X;2,0,0), e_{5}^3(X;1,1,0),$ $e_{5}^3(X;2,1,0), e_{5}^3(X;2,2,0), e_{5}^3(X;1,1,1), e_{5}^3(X;2,1,1), e_{5}^3(X;2,2,1),$ $e_{5}^3(X;2,2,2)\}$.}
\label{fig:Figure_11}
\end{figure}

By simple hand drawing with pencil and paper, we can play just for fun at this recreational math game, and so we did (for a couple of hours) streaming live on YouTube in January 2023.
The outcome was as follows: $d_{1}^4(Q)=0$, $d_{2}^4(Q)=1$, $d_{3}^4(Q)=2$, $d_{4}^4(Q)=2$, and $d_{5}^4(Q)=3$, but it is not easy going forward without a computer analysis (e.g., in order to find that $d_{5}^3(Q)=2$ we need to check, at least, $10$ different starting positions, as ${\rm{V_1}} \equiv (0,0,0)$, ${\rm{V_2}} \equiv (1,0,0)$, ${\rm{V_3}} \equiv (2,0,0)$, ${\rm{V_4}} \equiv (1,1,0)$, ${\rm{V_5}} \equiv (2,1,0)$, ${\rm{V_6}} \equiv (2,2,0)$, ${\rm{V_7}} \equiv (1,1,1)$, ${\rm{V_8}} \equiv (2,1,1)$, ${\rm{V_9}} \equiv (2,2,1)$, ${\rm{V_{10}}} \equiv (2,2,2)$, see Figure \ref{fig:Figure_11}), so we prefer to let the reader choose the most effective strategy to prove or disprove the general conjecture introduced in Section \ref{sec:2} stating that, for every given $k \in \mathbb{N}-\{0,1\}$, $\exists n^{\star}:=n^{\star}(k) : d_{n^{\star}+q}^k(Q)=d_{n^{\star}}^k(Q)=k \hspace{2mm} \forall q \in \mathbb{Z}^+$.

For $k=3$, the total number of starting positions that we have to check when we decide to find by brute force the value of $d_{n}^k(X)$ is $\binom{\lfloor \frac{n-1}{2}\rfloor +3}{3}$, which is equal to $\frac{2 \cdot n^3 + 15 \cdot n^2 + 31 \cdot n + 15 + 3 \cdot (-1)^{n+1} \cdot (n^2 + 5 \cdot n + 5)}{96}$ (i.e., the $(n-1)$-th duplicated tetrahedral number \cite{14}).

In general, we have that the number of starting vertices to be checked by brute force is the $(n-1)$-th duplicated $k$-simplicial number, and this value is equal to $\binom{\lfloor \frac{n-1}{2}\rfloor + k}{k}$ (the $k$-simplicial number sequence, $a_{k}(n)$, is defined as $a_{k}(n):=\frac{n \cdot (n+1)\cdot (n+2)\cdots (n+k-1)}{k!}$ since the $n$-th triangular number \cite{31} corresponds to the sum of the smallest $n$ natural numbers, the $n$-th tetrahedral number is the sum of the first $n$ triangular numbers \cite{32}, and so forth \cite{33}).
As an example, if $k=4$, then the number of cases that we need to check is equal to $\binom{\lfloor \frac{n-1}{2}\rfloor + 4}{4}$, which can be also written as $\frac{1}{24} \cdot (\floor{\frac{n - 1}{2}} + 1) \cdot (\floor{\frac{n - 1}{2}} + 2) \cdot (\floor{\frac{n - 1}{2}} + 3) \cdot (\floor{\frac{n - 1}{2}} + 4)$.

Basically, even if $k=3$ is given, the consideration above makes the complexity of the described raw approach harder and harder as $n$ grows, but it is still doable as $X=Q$ since $d_{6}^3(Q)=d_{6}^3(R)=3$ (in order to provide a direct proof that $d_{6}^3(Q) = 3$, we can take a look at the distance between the vertex $(0,0,0)$ and $(4,5,2) \in C(6,3)$, observing at the same time that $d_{6}^3(R) \leq 3$ holds by definition).

More in detail, we have that $e_{1}^k(X;0)=0 \Rightarrow r_{1}^k(X)=d_{1}^k(X)=0$ for any given chess piece $X$ among our canonical set of six.

\sloppy Now, $e_{2}^k(Q;0,0, \dots, 0)=1 \Rightarrow r_{2}^k(Q)=d_{2}^k(Q)=1$ and $(e_{3}^k(Q;0,0, \dots, 0, 0)=e_{3}^k(Q;1,0, \dots, 0, 0)=e_{3}^k(Q;1,1, \cdots, 0, 0)= \dots =e_{3}^k(Q;1,1, \dots, 1,0)=2 \wedge e_{3}^k(Q;1,1, \dots, 1,1)=1) \Rightarrow \bigl(r_{3}^k(Q)=1 \wedge d_{3}^k(Q)=2 \bigr)$ are also trivial (these relations hold by construction for any $k \geq 2$).

Then, assume $k=3$ and let us check what happens as $n : n > 3$ grows up to $n^{\star}$ (for every given $k \geq 2$, our conjecture implies that $d_{n^{\star}}^k(Q)=d_{n}^k(R)$ as long as $n \geq n^{\star}$).

\begin{figure}[H]
\begin{center}
\includegraphics[width=\linewidth]{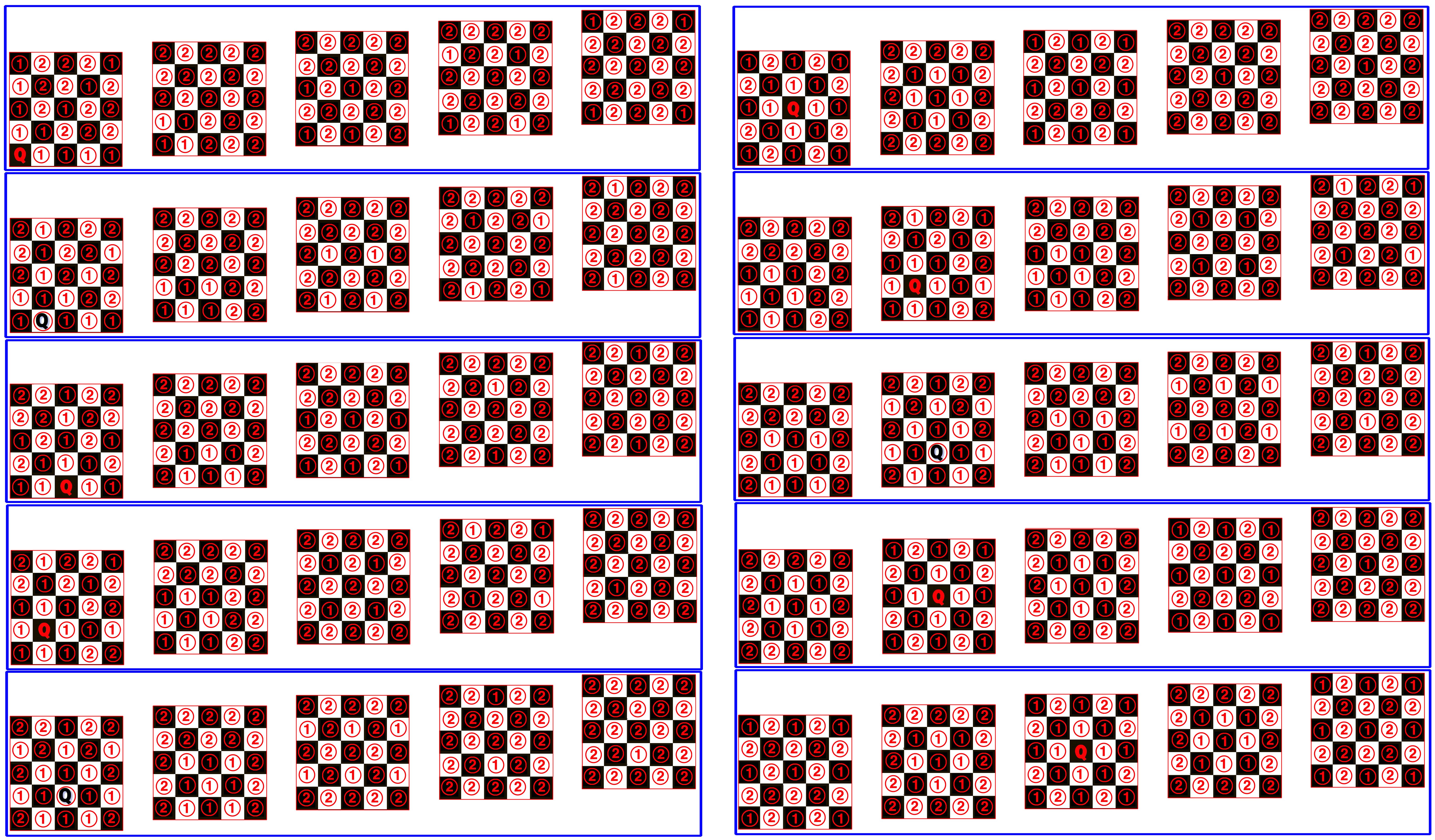}
\end{center}
\caption{Constructive proof that $r_{5}^3(Q)=2$ and $d_{5}^3(Q)=2$ by checking all the $10$ configurations represented in Figure \ref{fig:Figure_11}.}
\label{fig:Figure_12}
\end{figure}

From $\delta_{3}^3(R)((0,0,0),(2,2,1))=3$, since $e_{5}^3(Q;0,0,0)=e_{5}^3(Q;1,0,0)=e_{5}^3(Q;2,0,0)=e_{5}^3(Q;1,1,0)=e_{5}^3(Q;2,1,0)=e_{5}^3(Q;2,2,0)=e_{5}^3(Q;1,1,1)=e_{5}^3(Q;2,1,1)=e_{5}^3(Q;2,2,1)=e_{5}^3(Q;2,2,2)=2 \Rightarrow r_{5}^3(Q)=d_{5}^3(Q)=2$ (see Figure \ref{fig:Figure_12}), and by observing that the queen distance between $(0,0,0)$ and $(4,5,2) \in C(6,3)$ is $3$, $\min\{n :d_{n}^3(Q)=3\}=6$ follows (i.e., if $k=3$ is given, then $n^{\star}=6$).

\begin{figure}[H]
\begin{center}
\includegraphics[width=\linewidth]{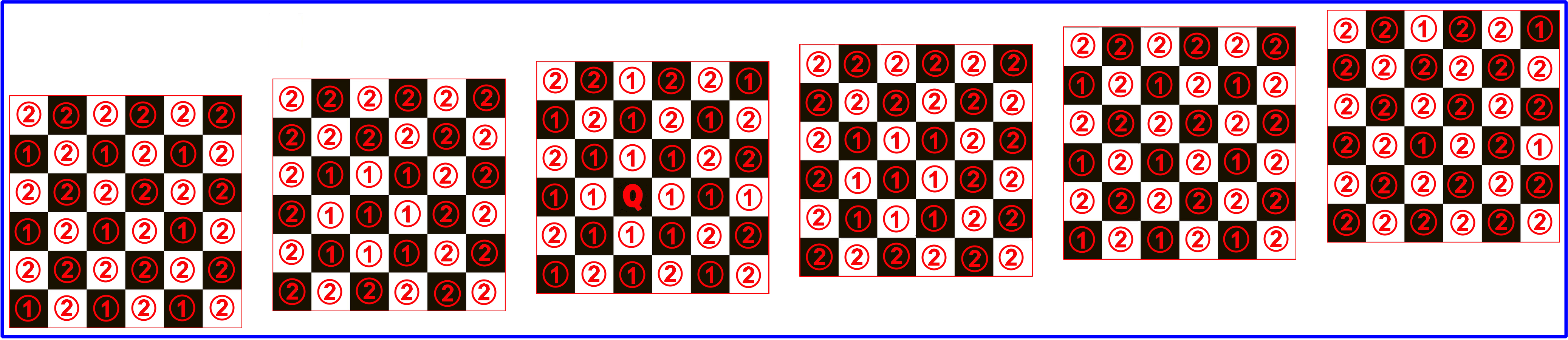}
\end{center}
\caption{Graphical proof that $r_{6}^3(Q)=2$ since $\nexists (x_1,x_2,x_3) : S_1(Q;6,3;x_1,x_2,x_3) = C(6,3)$.}
\label{fig:Figure_13}
\end{figure}

\begin{figure}[H]
\begin{center}
\includegraphics[width=\linewidth]{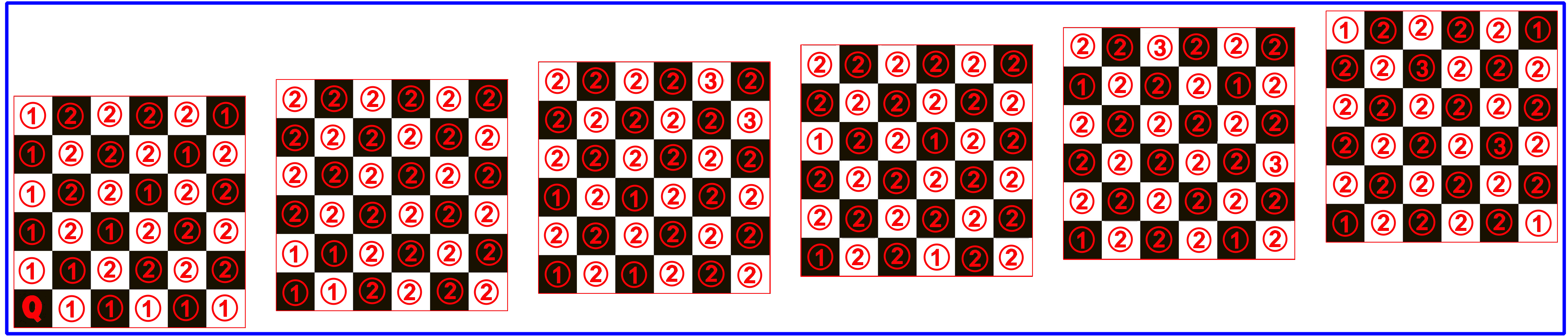}
\end{center}
\caption{Graphical proof that $n \geq 6 \Rightarrow d_{n}^3(Q)=3$ since $d_{n}^k(Q) \leq k$ by construction.}
\label{fig:Figure_14}
\end{figure}

Although from $e_{6}^3(Q;2,2,2)=2$ (see Figure \ref{fig:Figure_13}) we deduce that $r_{6}^3(Q)=2$, and thus \linebreak $r_{6}^3(Q)<d_{6}^3(Q)$, we have finally concluded that $d_{n}^3(Q)=3$ for any $n \geq 6$ (see Figure \ref{fig:Figure_14}).

The next step is to study $d_{n}^4(Q)$, bearing in mind the upper bound $d_{n}^4(Q) \leq 4$. By repeating the previous (tedious) approach for $n=4$, we have found that $d_{4}^4(Q)=r_{4}^4(Q)=2$, but this is not a valid method to determine the general formula that describes the sequence \linebreak $a(k):=\min\{n \in \mathbb{Z}^+ :d_{n}^k(Q)=k\}$, for all $k \in \mathbb{N}-\{0,1\}$. Thus, the aforementioned problem is entirely open and suitable for recreational mathematicians, computer scientists, graph theorists, and many others.

The pawn is a very special chess piece since the well-known Article 3.7, letter ${\rm{e}}$, states that it can promote to another piece belonging to the set $\{Q,R,N,B\}$, but for our purposes we can reduce it to just $\{Q,N\}$, given the fact that our final goal is to evaluate the smallest pawn distance between any pair of vertices of $C(n,k)$.

In the present section, we consider only the minimum number of pawn moves necessary to go from a given vertex of $C(n,k)$ to another one belonging to the same set. Now, since $\bigl(S_t(B;n,k;x_1, x_2, \dots, x_k) \subseteq S_t(Q;n,k;x_1, x_2, \dots, x_k)\bigr) \wedge \bigl(S_t(R;n,k;x_1, x_2, \dots, x_k) \subseteq S_t(Q;n,k;x_1, x_2, \dots, x_k)\bigr)$ is true for any $n-1, \hspace{1.5mm} k-2, \hspace{1.5mm} t \hspace{-2.5mm} \in \hspace{-3mm} \mathbb{N}_0$, it follows that $d_{n}^k(Q) \leq \min\{d_{n}^k(R),d_{n}^k(B)\}$ and also $r_{n}^k(Q) \leq \min\{r_{n}^k(R),r_{n}^k(B)\}$ so that the only open match comes from promoting to a queen or a knight. This outcome depends on the fact that $\bigl(S_t(N;n,k;x_1, x_2, \dots, x_k) \nsubseteq S_t(Q;n,k;x_1, x_2, \dots, x_k)\bigr) \wedge \bigl(S_t(N;n,k;x_1, x_2, \dots, x_k) \nsupseteq S_t(Q;n,k;x_1, x_2, \dots, x_k)\bigr)$ and so we cannot lightly state that, in general, \linebreak$e_{n}^k(Q;x_1, n-1, \dots, x_k) \leq e_{n}^k(N;x_1, n-1, \dots, x_k)$ (or vice versa).

Given $X:=P$, we would always expect to observe that promoting to a queen cannot be worse than opting for a knight promotion. On the other hand, we cannot infer that, in chess, the minimum distance between any vertex ${\rm{V_1}} \in C(n,k)$ and another vertex ${C(n,k) \ni \rm{V_2}} \not\equiv {\rm{V_1}}$ is a queen distance, instead of a knight distance.

By assuming $n=3$ and $k=2$, Figure \ref{fig:Figure_15} proves that a pawn can reach any square on the board in $2$ moves or less, improving the general upper bound $d_{n}^k(Q)+n-2$ for the given case (since $d_{3}^2(Q)+3-2=2+3-2=3$), but this surprising outcome can be achieved if and only if the pawn is free to alternatively promote to a queen or knight, depending on which path is better for the given configuration.

\begin{figure}[H]
\begin{center}
\includegraphics[width=\linewidth]{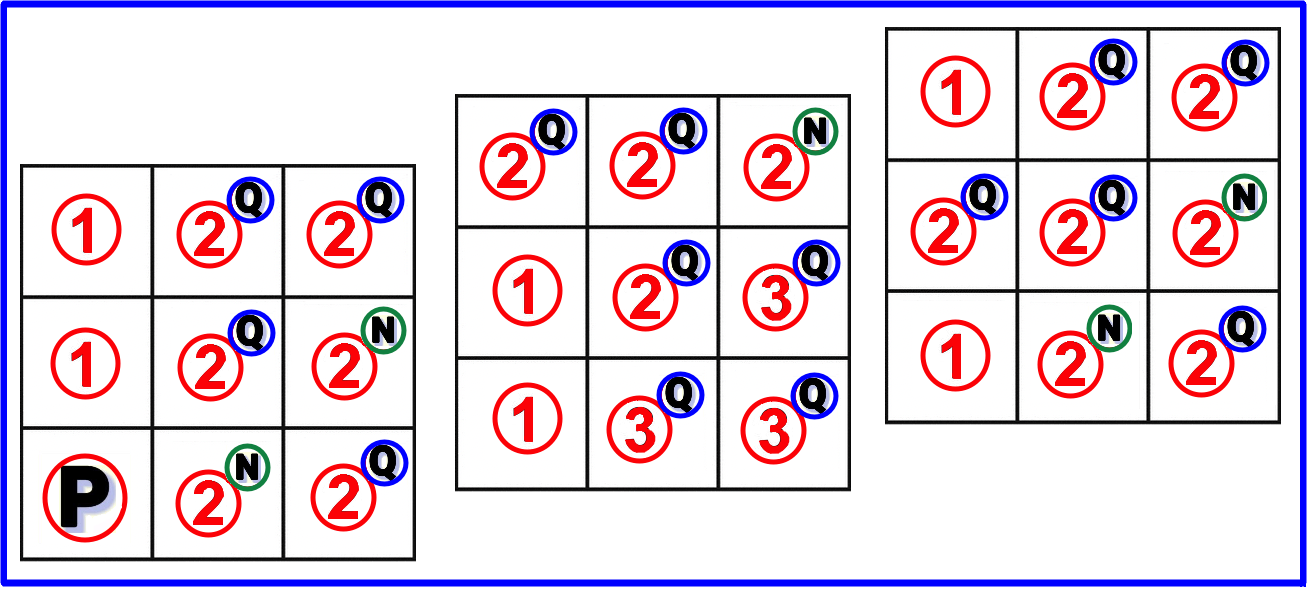}
\end{center}
\caption{Graphical proof that $2$ moves are not sufficient to let a pawn reach any square of the $3 \times 3 \times 3$ chessboard when we set $(0,0,0)$ as the starting vertex, whereas $3$ moves are enough for the same board since the trivial upper bound $d_{n}^k(Q)+n-2$ cannot be violated.}
\label{fig:Figure_16}
\end{figure}

As shown by Figure \ref{fig:Figure_16}, we observe that the above-mentioned trick does not work anymore if $k$ is greater than $2$.

\begin{figure}[H]
\begin{center}
\includegraphics[width=\linewidth]{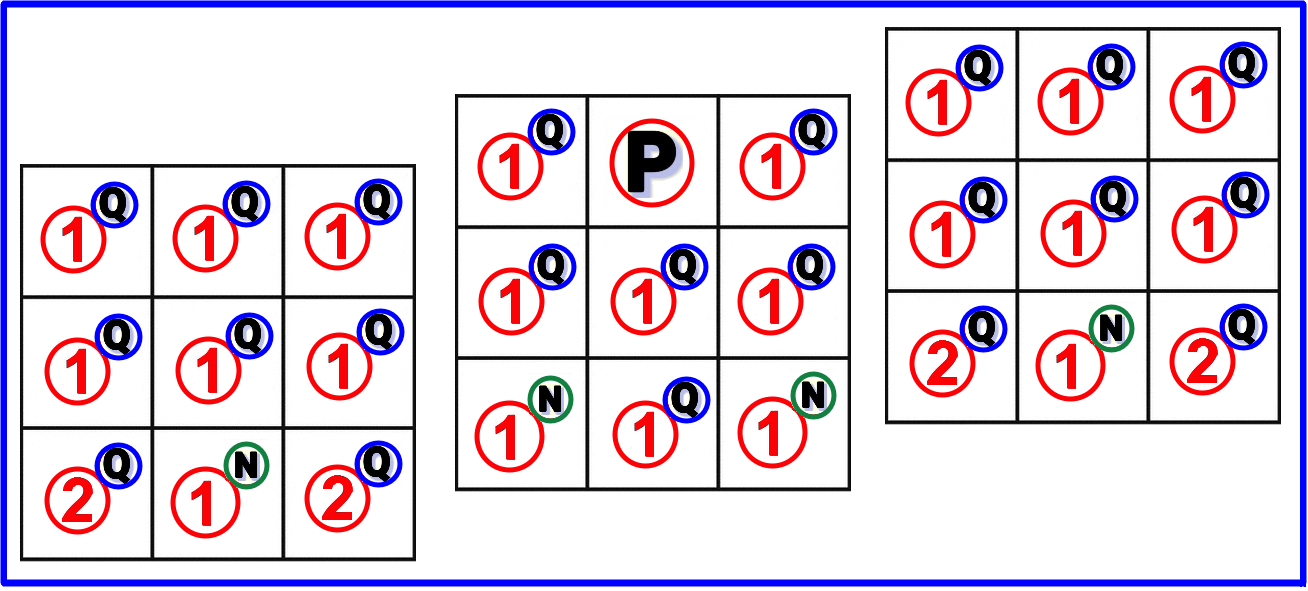}
\end{center}
\caption{Number of required moves for a pawn (i.e., $P$ or $\overline{P}$) placed on the vertex $(1,2,1)$ to reach any other vertex of $C(3,3)$.}
\label{fig:Figure_0}
\end{figure}

Moreover, by taking the point $(1,2,1)$ as our starting vertex, it is not difficult to check that $2$ moves are sufficient to let the pawn reach any other vertex of $C(3,3)$ if we admit, as we did, that starting from a vertex belonging to the last rank let us promote (to a queen or knight) at move zero. This consideration shows that a pawn can always reach any other vertex of $C(n,k)$ in no more than $d_{n}^k(Q)$ moves, as long as the mentioned pawn promotes to a queen and is left free to start from the $(n-1)$-th rank (see Figure \ref{fig:Figure_0}).

In conclusion, as $n$ or $k$ grows, we need to carefully check case by case \linebreak$S_t(N;n,k;x_1, n-1, x_3, \dots, x_k)$ and $S_t(Q;n,k;x_1, n-1, x_3, \dots, x_k)$ (up to $t=d_{n}^k(Q)$, if necessary) in order to confidently state that $d_{n}^k(Q)+n-2$ is the minimum number of $\overline{P}$-moves needed to reach the farthest vertex of $C(n,k)$ from ${C(n,k) \ni \rm{V_1}} \equiv (x_1,x_2,x_3,\dots,x_k)$, by assuming the worst case scenario.
 
Lastly, it is time to study the sixth and last chess piece of our set, the fascinating $k$-dimensional knight (a.k.a. $k$-knight).

In this case, for any $k \geq 2$, we know in advance that $d_{1}^k(N)=r_{1}^k(N)=0$, $n \in \{2,3\} \Rightarrow r_{n}^k(N)=+\infty$, and
$n \geq 4 \Rightarrow d_{n}^3(N)=n$ (while $d_{n}^2(N)$ is fully described by the OEIS sequence A232007 \cite{15}). So, from here on, let us assume that $n,k \in \mathbb{N}-\{0,1,2,3\}$.

About the gap $d_{n}^3(N) - e_{n}^3(N;x_1, x_2, x_3)$, it depends on $(x_1, x_2, x_3)$ itself.
As a couple of examples, we have that $4=e_{4}^3(N;1, 0, 0) \wedge 3=e_{4}^3(N;1, 1, 0)=e_{4}^3(N;1, 1, 1)$ implies $d_{4}^3(N)-r_{4}^3(N) = 1$, while $5=e_{5}^3(N;0, 0, 0)=e_{5}^3(N;1, 0, 0) \wedge 4=e_{5}^3(N;2, 0, 0)=e_{5}^3(N;1, 1, 0)=e_{5}^3(N;2, 1, 0)=e_{5}^3(N;2, 2, 0)=e_{5}^3(N;1, 1, 1)=e_{5}^3(N;2, 1, 1) \wedge 3=e_{5}^3(N;2, 2, 1)=e_{5}^3(N;2, 2, 2)$ implies that $d_{5}^3(N)-r_{5}^3(N) = 2$, and so forth (e.g., $e_{6}^3(N;0, 0, 0)=1+e_{6}^3(N;2, 2, 0)=2+e_{6}^3(N;2, 2, 2)$ so that the difference $d_{6}^3(N)-r_{6}^3(N)$ cannot be smaller than $d_{5}^3(N)-r_{5}^3(N)$).

When $k$ is an arbitrarily large integer, we invoke the lower bound arising from Thomas Andrews' answer to the MSE question \#$4619708$ \cite{18}, which shows that $e_{n}^k(N;0, 0, \dots, 0) \geq \lfloor \frac{k \cdot (n-1)}{3} + 1 \rfloor$ generally holds and, if $k \cdot (n-1)$ belongs to the congruence class $2$ modulo $3$, then the aforementioned bound can be improved as $e_{n}^k(N;0, 0, \dots, 0) \geq \lfloor \frac{k \cdot (n-1)}{3} + 2 \rfloor$. It is worth mentioning the fact that these are just lower bounds for $e_{n}^k(N;0, 0, \dots, 0)$ since Andrews' answer did not take into account the constraint, stated in the present paper, that the knight can only move between elements of the set $C(n,k)$ so that no (Steiner) point $(x_1, x_2, \dots, x_k) : \linebreak (x_1 \notin \{0,1,\dots,n-1 \}) \vee (x_2 \notin \{0,1,\dots,n-1 \}) \vee \dots \vee (x_k \notin \{0,1,\dots,n-1 \})$ can be considered for any knight path in order to study $d_{n}^k(N)$ or $r_{n}^k(N)$.
As a clarifying example, let $n \geq 4$ be given and observe that $(0,0)$-$(2,-1)$-$(1,1)$ is not a valid path (since the $x_2$ coordinate of $(2,-1)$ is a negative number), whereas a planar minimal path from $(0,0)$ to $(1,1)$ is achievable as $(0,0)$-$(2,1)$-$(1,3)$-$(3,2)$-$(1,1)$, doubling the knight distance from $(0,0,0)$ to $(1,1,0)$ (since $(0,0,0)$-$(1,0,2)$-$(1,1,0)$ describes a valid $3$D path in $[0,1] \times [0,1] \times [0,2]$).

Thus, we can state that $d_{n}^k(N) \geq \lfloor \frac{k \cdot (n-1)}{3} + 1 \rfloor$ 
holds for any $n,k \in \mathbb{N}-\{0,1,2,3\}$, and the comparison between the lower bound $d_{n}^k(N) \geq \lfloor \frac{k \cdot (n-1)}{3} + 1 \rfloor$ and the raw upper bound $d_{n}^k(Q) \leq k$ explains the main reason behind our conjecture that $(n \neq 3 \vee k \neq 2)$ represents a sufficient condition for requiring exactly $d_{n}^k(Q)+n-2$ pawn moves to reach the farthest vertex of $C(n,k)$, by assuming a worst-case configuration.


\section{Defining the true international chess pieces in higher dimensions} \label{sec:4}

In Subsection \ref{sec:sub4.1} we suggest how to extend the definition of every ($k$-dimensional) international chess piece in the most strict and coherent mathematical way, despite the (pretty arbitrary) choices we made in the previous sections by applying retrograde check to the complexity of the possible outcomes.

Then, in Subsection \ref{sec:sub4.2}, we will study the radius and the diameter of each of the graphs induced by the six generalized $k$-dimensional chess pieces stated in the next subsection.

\subsection{Coloring any multidimensional chessboard from a mathematical perspective}
\label{sec:sub4.1}

To properly define the $k$-dimensional versions of the six FIDE chess pieces, let us consider the issue related to coloring the set $C(n,k)$ in a way that a generalized light-square bishop (a.k.a. odd $k$-bishop) can reach all and only the vertices of $C(n,k)$ colored as its own starting vertex since this is a key rule of traditional chess that every player needs to learn before starting his very first game of chess.

Furthermore, we also need to take into account the definition of $k$-queen, which is a piece that alternatively moves as a bishop or a rook \cite{1}, and the fact that if we put a rook on the vertex $(x_1, x_2, \dots, x_k) \in C(n,k)$ and perform a single move, it cannot overlap any vertex of $C(n,k)$ that a bishop placed on $(x_1, x_2, \dots, x_k)$ can reach in its next move.

For these reasons, let us indicate as $\tilde{B}$ our new $k$-bishop, as $\tilde{R}$ our improved $k$-rook, and as $\tilde{Q}$ the resulting version of $k$-queen. They must satisfy the following three conditions.

\begin{enumerate}
\item $S_0(\tilde{B};n,k;x_1, x_2, \dots, x_k)=S_0(\tilde{R};n,k;x_1, x_2, \dots, x_k)=S_0(\tilde{Q};n,k;x_1, x_2, \dots, x_k)=\{(x_1, x_2, \dots, x_k)\}$.
\item $S_1(\tilde{B};n,k;x_1, x_2, \dots, x_k) \cap S_1(\tilde{R};n,k;x_1, x_2, \dots, x_k)=S_0(\tilde{Q};n,k;x_1, x_2, \dots, x_k)$.
\item $S_1(\tilde{B};n,k;x_1, x_2, \dots, x_k) \cup S_1(\tilde{R};n,k;x_1, x_2, \dots, x_k)=S_1(\tilde{Q};n,k;x_1, x_2, \dots, x_k)$.
\vspace{3mm}
\end{enumerate}
\vspace{-5mm}

On the other hand, reading Article 3.4 \cite{1} of the mentioned FIDE Handbook in conjunction with Article 3.6, we observe that we may also include the statement
\begin{enumerate}
\setcounter{enumi}{3}
\vspace{-3mm}
\item $S_1(\tilde{Q};n,k;x_1, x_2, \dots, x_k) \cap S_1(\tilde{N};n,k;x_1, x_2, \dots, x_k)=\{(x_1, x_2, \dots, x_k)\}$ 
\end{enumerate}
\vspace{-3mm}
as an additional constraint.

Accordingly, let us formally describe the structure of the new battlefield so that we will be finally ready to define the $k$-bishop ($\tilde{B}$), the $k$-rook ($\tilde{R}$), and also the $k$-queen ($\tilde{Q}$). 

\begin{definition} \label{def5.1}
The six generalized chess pieces (i.e., the $k$-bishop $\tilde{B}$, the $k$-rook $\tilde{R}$, the $k$-queen $\tilde{Q}$, the $k$-king $\tilde{K}$, the $k$-knight $\tilde{N}$, and the $k$-pawn $\tilde{P}$) move on a colored set of $n^k$ points, \linebreak$C(n,k):=\{\{0,1,\dots,n-1\}\times\{0,1,\dots,n-1\}\times \cdots \times\{0,1,\dots,n-1\}\} \subseteq \mathbb{Z}^k$, where the subset of the dark vertices of $C(n,k)$ is given by the elements that are characterized by an even value of $x_1 + x_2 + \cdots + x_k$, whereas the light vertices are those such that $\sum_{j=1}^{k} x_j=2 \cdot m +1, m \in \mathbb{N}_0$.
\end{definition}

As a result, any vertex ${\rm{V}} \equiv (x_1, x_2, \dots, x_k)$ belonging to the set $C(n,k)$ is a light vertex or a dark vertex and, by Definition \ref{def5.1}, we have that ${\rm{V}}$ is a light vertex if and only if $\sum_{j=1}^{k} x_j=2 \cdot m +1, m \in \mathbb{N}_0$. Otherwise, ${\rm{V}} \equiv (x_1, x_2, \dots, x_k)$ is a dark vertex since $\sum_{j=1}^{k} x_j=2 \cdot m, m \in \mathbb{N}_0$.

\begin{figure}[H]
\begin{center}
\includegraphics[width=\linewidth]{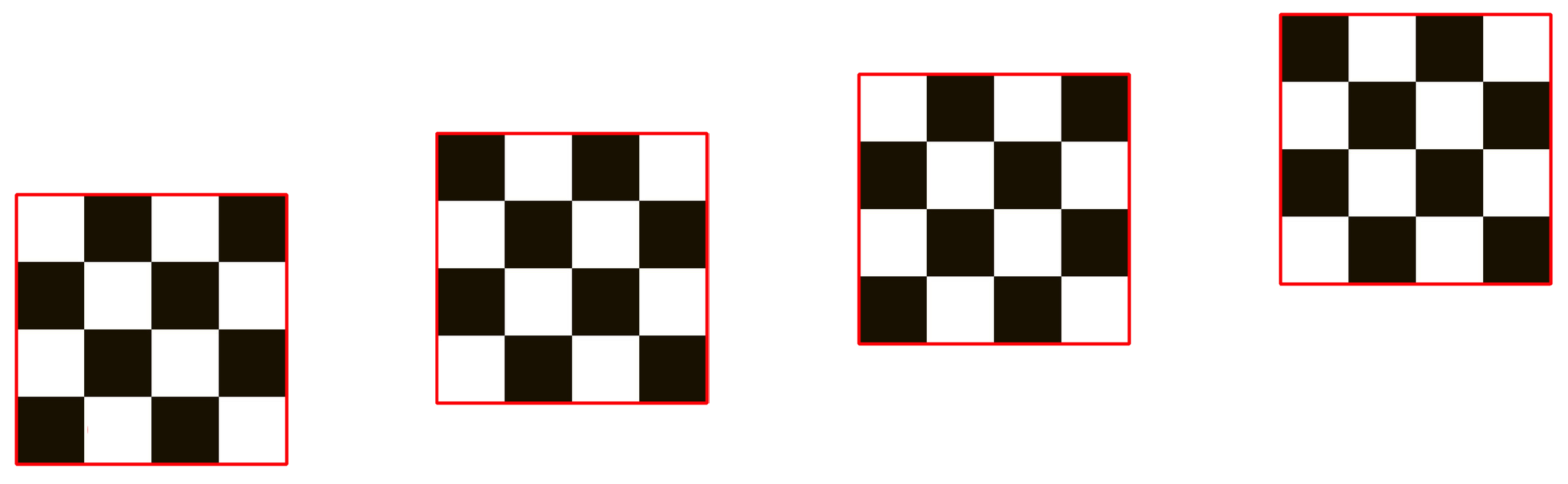}
\end{center}
\caption{Coloring the $4 \times 4 \times 4$ chessboard, consistently with the $k$-bishop ($\tilde{B}$) definition provided in the present section.}
\label{fig:Figure_17}
\end{figure}

This means that our $k$-dimensional checkerboard cannot be colored as shown in Figures \ref{fig:Figure_1} to \ref{fig:Figure_9}, but it looks as follows (Figure \ref{fig:Figure_17}).

We can see that this is the correct $C(n,k)$ setup by taking the previously stated knight definition (i.e., $S_1(N;n,k;x_1, x_2, \dots, x_k)=\{(x_1, x_2, \dots, x_k)\} \cup \{(x_1+c_1, x_2+c_2, \dots, x_k+c_k) : \forall j \in \{1,2,\dots,k\} \left( \sum_{j=1}^{k} |c_j| = 3 \wedge |c_j| \in \{0,1,2\} \wedge \exists! j : |c_j|=1 \right) \wedge (x_j+c_j) \in \{0,1,\dots,n-1 \}\}$, see Figure \ref{fig:Figure_7}) and checking that the $k$-knight always moves from a dark to a light square and vice versa, as it always occurs on the standard $8 \times 8$ chessboard when an official FIDE chess game is played \cite{1}.
Thus, let $\tilde{N}:=N$.

Now, we can properly define the $k$-bishop ($\tilde{B}$).
Let $x_1, x_2, \dots, x_k \in \{0,1,\dots,n-1\}$ be given. Then, $S_0(\tilde{B};n,k;x_1, x_2, \dots, x_k):=\{(x_1, x_2, \dots, x_k)\}$ and $S_1(\tilde{B};n,k;x_1, x_2, \dots, x_k)=\{(x_1, x_2, \dots, x_k)\} \cup \{(x_1+c_1, x_2+c_2, \dots, x_k+c_k) : \exists c \in \{1,\dots,n-1\}\quad \forall j \in \{1,2, \dots, k\} \quad c_j \in \{-c, 0, c \} \wedge |\{j : c_j \neq 0 \}|=2 \wedge (x_j+c_j) \in \{0,1,\dots,n-1 \}\}$ (i.e., $\{S_1(\tilde{B};n,k;x_1, x_2, \dots, x_k)-S_0(\tilde{B};n,k;x_1, x_2, \dots, x_k)\}$ includes all the vertices $(x_1+c_1, x_2+c_2, \dots, x_k+c_k)$ of $C(n,k)$ such that exactly two of the coordinate increments $c_1,c_2,\dots,c_k$ are nonzero).

\begin{figure}[H]
\begin{center}
\includegraphics[width=\linewidth]{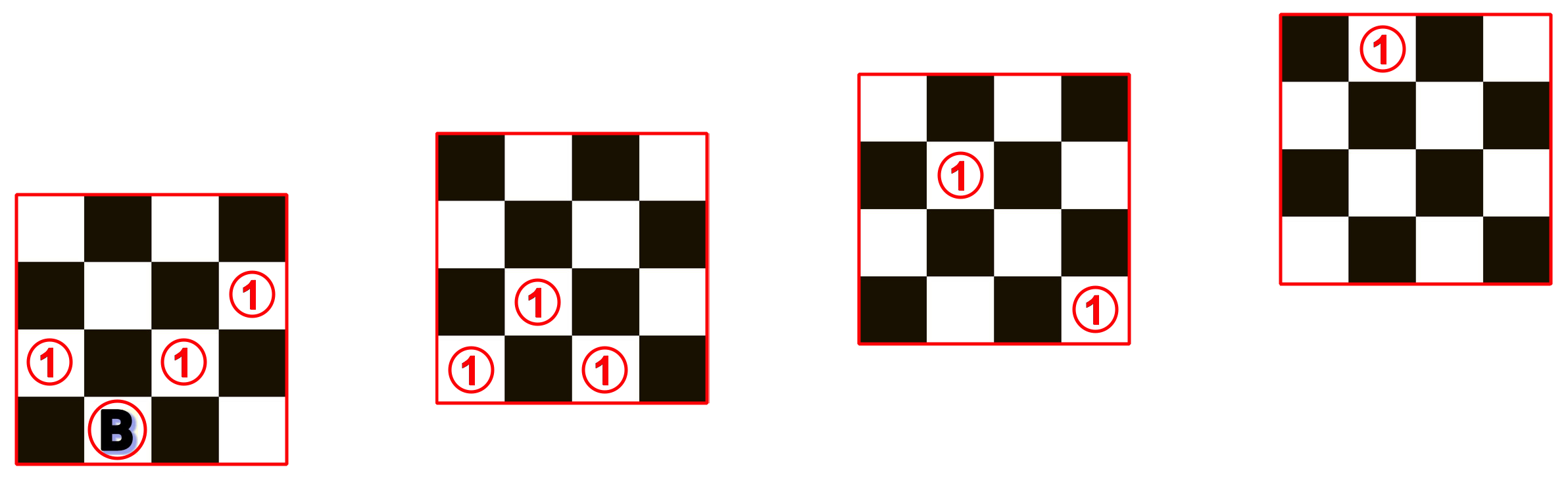}
\end{center}
\caption{The set $S_1(\tilde{B};n,k;x_1, x_2, \dots, x_k)$ represents the most natural way to generalize a bishop in $\mathbb{Z}^k$ since it moves only on the vertices such that the sum of their $k$ coordinates is always odd if $\sum_{j=1}^{k} x_j$ is odd (and vice versa if $\sum_{j=1}^{k} x_j$ is even).}
\label{fig:Figure_22}
\end{figure}

Figure \ref{fig:Figure_22} shows the set $S_1(\tilde{B};4,3;1,0,0) = \{(1,0,0), (0,1,0), (2,1,0), (3,2,0), (0,0,1),$ $(1,1,1), (2,0,1), (1,2,2), (3,0,2), (1,3,3)\}$ originated by choosing the vertex $(1,0,0) \in C(4,3)$ as our starting point (i.e., $S_0(\tilde{B};4,3;1,0,0)=\{(1,0,0)\}$).

\begin{figure}[H]
\begin{center}
\includegraphics[width=\linewidth]{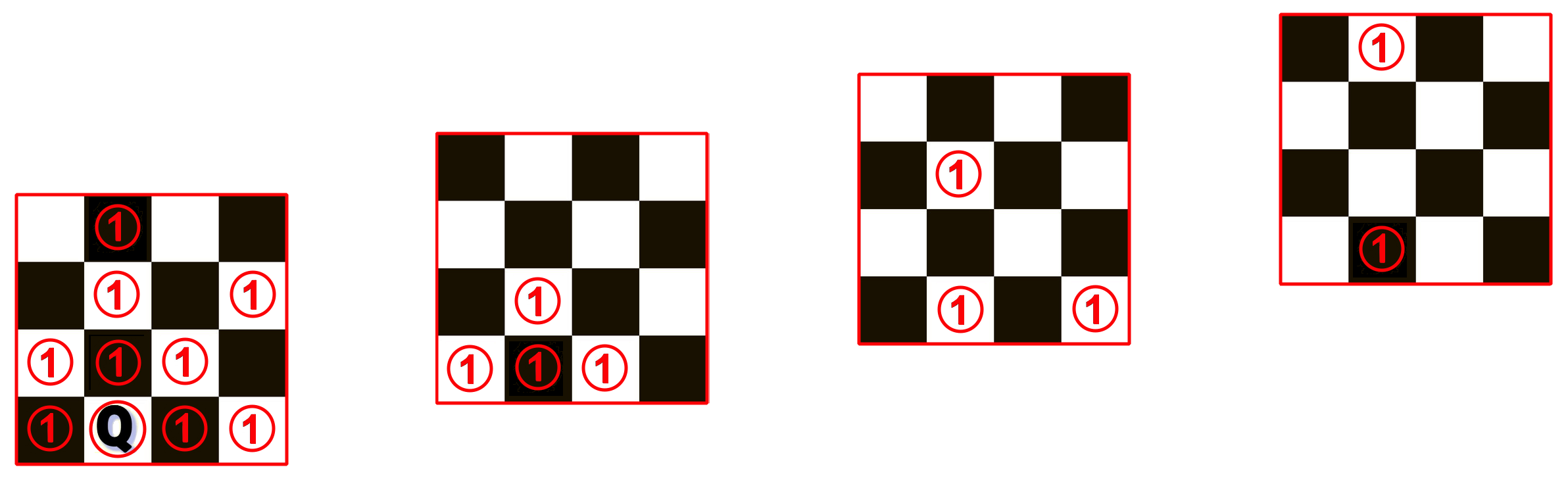}
\end{center}
\caption{The set $S_1(\tilde{Q};n,k;x_1, x_2, \dots, x_k)$ is the right way to generalize a queen in $\mathbb{Z}^k$ since it incorporates the moving pattern of a $k$-bishop and a $k$-rook, with only the starting vertex in the intersection of the two sets, $S_1(\tilde{B};n,k;x_1, x_2, \dots, x_k)$ and $S_1(\tilde{R};n,k;x_1, x_2, \dots, x_k)$, as suggested by the third article of the FIDE Handbook. Furthermore, we observe that the condition $S_1(\tilde{Q};n,k;x_1, x_2, \dots, x_k) \cap S_1(\tilde{N};n,k;x_1, x_2, \dots, x_k) = \{(x_1, x_2, \dots, x_k)\}$ is also satisfied.}
\label{fig:Figure_30}
\end{figure}

By looking at the description of $3$-queen provided by the mentioned Millennium $3$D chess (and taking into account how its $3$D board has been colored there), we cannot agree anymore that Figure \ref{fig:Figure_5} provides a coherent extension of the (standard) FIDE chess queen to the next dimension. Instead, Figure \ref{fig:Figure_30} shows the most clean and rational way to generalize a FIDE $2$-queen to higher dimensions.

As a consequence, we cannot confirm here the generalization of the FIDE queen suggested in Section \ref{sec:Intr} of the present paper (as it was derived by assuming that $S_1(B;n,k;x_1, x_2, \dots, x_k)=S_0(B;n,k;x_1, x_2, \dots, x_k) \cup \{(x_1+c_1, x_2+c_2, \dots, x_k+c_k) : \forall j \in \{1,2,\dots,k\} \linebreak (c_j \in \{-c,0,c\} \wedge c_j \neq 0 \hspace{2mm} \textnormal{if} \hspace{2mm}j < 3) \wedge (x_j+c_j) \in \{0,1,\dots,n-1\}\}$ describes the moving rule of any multidimensional bishop), since we have that $S_1(\tilde{Q};n,k;x_1,x_2,\dots,x_k):=S_1(\tilde{B};n,k;x_1,x_2,\dots,x_k) \cup S_1(\tilde{R};n,k;x_1,x_2,\dots,x_k)$ and so $S_1(\tilde{Q};n,k;x_1,x_2,\dots,x_k) = S_1(\tilde{B};n,k;x_1, x_2, \dots, x_k) \cup S_1(\overline{R};n,k;x_1, x_2, \dots, x_k)$, where $S_1(\overline{R};n,k;x_1,x_2,\dots,x_k)$ is shown in Figure \ref{fig:Figure_4}. Thus, for every positive integer $t$, we have that $S_t(\tilde{R};n,k;x_1, x_2, \dots, x_k) = S_t(\overline{R};n,k;x_1,x_2,\dots,x_k)$.

From $\tilde{R}:=\overline{R}$ (since $\overline{R}$, by definition, changes one and only one coordinate at each move) we note that the $k$-rook moves from a light vertex to a dark vertex of $C(n,k)$ if and only if $\sum_{j=1}^{k} x_j$ and $c : 1 \leq |c| \leq n-1$ are both odd (see Figure \ref{fig:Figure_17}). This is perfectly consistent with the well-known moving pattern of the FIDE planar rook \cite{1}.

Moreover, for any $n,k \in \mathbb{N}-\{0,1\}$, the equivalence $d_{n}^k(\tilde{R})=r_{n}^k(\tilde{R})=k$ trivially follows from the $k$-rook definition itself.

Since the standard FIDE pawn is not allowed to move backward, it cannot induce any metric space as it is (i.e., even the second condition stated at the beginning of Section \ref{sec:Intr} would be violated). So, in order to change this ugly chess asymmetry, we choose to remove the aforementioned limitation by defining the $k$-pawn distance between any pair of distinct vertices of $C(n,k)$, ${\rm{V_1}} \equiv (x_1,x_2,x_3,\ldots,x_k)$ and ${\rm{V_2}} \equiv (x_1+c_1,x_2+c_2,x_3+c_3,\ldots,x_k+c_k)$ such that $x_1, x_1+c_1 ,x_2, x_2+c_2, x_3, x_3+c_3, \ldots,x_k, x_k+c_k \in \{0,1,2,\dots,n-1\}$, as $\delta_{n}^k(\tilde{P})({\rm{V_1}}, {\rm{V_2}}) := n-2+\max\{d_{n}^2(\tilde{Q}), d_{n}^2(\tilde{N}) \}+\min\{\delta_{n}^k(\tilde{Q})({\rm{V_1}}, {\rm{V_2}}),\delta_{n}^k(\tilde{N})({\rm{V_1}}, {\rm{V_2}}) \}$, for any $n \in \mathbb{N}-\{0,1,2,3\}$ and $k \in \mathbb{N}-\{0,1\}$.

The key idea here is to consider a standard FIDE pawn which lies on the first rank of a $k$-dimensional board, in generic position $(x_1,0,x_3,\ldots,x_k)$. 

Now, according to Article 3.7 \cite{1}, any pawn can step forward by two squares while performing its first move and then it moves forward by one square only.

Thus, in the worst-case scenario, our $k$-pawn is placed on the first rank, in position $(x_1,0,x_3,\ldots,x_k)$, and it needs $n-2$ moves to reach the last rank (i.e., the vertex \linebreak $(x_1,n-1,x_3,\ldots,x_k)$) and promote to a queen or a knight. Then, since we are assuming $n > 3$, the new queen would spend one additional move to come back to $(x_1,0,x_3,\ldots,x_k)$ while the $k$-knight will certainly be able to do the same in (at most) $d_{n}^2(\tilde{N})$ additional moves. On this purpose, it is worth to note that, given $C(n,2)$ such that $n \geq 5$, a knight can always go from ${\rm{V_1}} \in C(n,2)$ to ${\rm{V_2}} \in C(n,2)$ in $\left\lceil{\frac{2 \cdot n}{3}}\right\rceil$ moves or less (since $n \geq 5 \Rightarrow d_{n}^2(\tilde{N})=\left\lceil{\frac{2 \cdot n}{3}}\right\rceil$, see the OEIS sequences \cite{15, 20}).

We complete the definition of the $k$-pawn distance in $C(n,k)$ by covering the remaining $3$ cases, so let $n \in \{1,2,3\}$ be given. For every $n$, we set $\delta_n^k(\tilde{P})({\rm{V}},{\rm{V}}):=0$. Then, the following rules fully describe the $k$-pawn distance, $\delta_{n}^k(\tilde{P})$, between ${\rm{V_1}} \in C(n,k)$ and ${\rm{V_2}} \in C(n,k)$.

$n=1 \Rightarrow {\rm{V_1}} \equiv {\rm{V_2}} \Rightarrow \delta_{1}^k(\tilde{P})({\rm{V_1}}, {\rm{V_2}})=\delta_{1}^k(\tilde{P})({\rm{V_1}}, {\rm{V_1}}):=0$;

$n=2 \wedge {\rm{V_1}}\ne {\rm{V_2}} \Rightarrow \delta_{2}^k(\tilde{P})({\rm{V_1}}, {\rm{V_2}}):=2-2+d_{2}^2(\tilde{Q})+\delta_{2}^k(\tilde{Q})({\rm{V_1}}, {\rm{V_2}})=1+\delta_{2}^k(\tilde{Q})({\rm{V_1}}, {\rm{V_2}})$, since $d_{2}^2(\tilde{Q})=1$.

$n=3 \wedge {\rm{V_1}}\ne {\rm{V_2}} \Rightarrow \delta_{3}^k(\tilde{P})({\rm{V_1}}, {\rm{V_2}}):=3-2+d_{3}^2(\tilde{Q})+\min\{\delta_{3}^k(\tilde{Q})({\rm{V_1}}, {\rm{V_2}}),\delta_{3}^k(\tilde{N})({\rm{V_1}}, {\rm{V_2}})\} = 3+\min\{\delta_{3}^k(\tilde{Q})({\rm{V_1}}, {\rm{V_2}}),\delta_{3}^k(\tilde{N})({\rm{V_1}}, {\rm{V_2}})\}$ since $d_{3}^2(\tilde{Q})=2$ (e.g., $\delta_{3}^2(\tilde{N})((0,0), (1,2))=1$, while $\delta_{3}^2(\tilde{N})((0,0), (1,1))=+\infty$, and thus it follows that $\left({\rm{V_1}}=(0,0)\wedge {\rm{V_2}}=(1,2) \right) \Rightarrow \delta_{3}^2(\tilde{P})({\rm{V_1}}, {\rm{V_2}})=3+\delta_{3}^2(\tilde{N})({\rm{V_1}}, {\rm{V_2}})$, whereas $\left({\rm{V_1}}=(0,0)\wedge {\rm{V_2}}=(1,1) \right) \Rightarrow \delta_{3}^2(\tilde{P})({\rm{V_1}}, {\rm{V_2}})=3+\delta_{3}^2(\tilde{Q})({\rm{V_1}}, {\rm{V_2}})$).

As a peculiar case, we can also observe how $n=4$ implies that $\delta_{4}^k(\tilde{P})({\rm{V_1}}, {\rm{V_2}})=4-2+d_{4}^2(\tilde{N})+\min\{\delta_{4}^k(\tilde{Q})({\rm{V_1}}, {\rm{V_2}}),\delta_{4}^k(\tilde{N})({\rm{V_1}}, {\rm{V_2}})\}$.

Thus, $\delta_{4}^k(\tilde{P})({\rm{V_1}}, {\rm{V_2}}) = 7+\min\{\delta_{4}^k(\tilde{Q})({\rm{V_1}}, {\rm{V_2}}),\delta_{4}^k(\tilde{N})({\rm{V_1}}, {\rm{V_2}}) \}
\leq 7+\delta_{4}^k(\tilde{Q})({\rm{V_1}}, {\rm{V_2}}) \leq 7+\delta_{4}^k(\tilde{R})({\rm{V_1}}, {\rm{V_2}})$.

Hence, $\delta_{4}^k(\tilde{P})({\rm{V_1}}, {\rm{V_2}}) \leq 7+k$ (since $\delta_{4}^k(\tilde{R})({\rm{V_1}}, {\rm{V_2}})\leq k$, for any $k > 1$).

In conclusion, let us state that if $n \geq 5$, $C(n,k) \ni {\rm{V_1}} \equiv (x_1, x_2, \ldots, x_k)$, $C(n,k) \ni {\rm{V_2}} \equiv (x_1+c_1,x_2+c_2,\ldots,x_k+c_k)$, and ${\rm{V_1}} \not\equiv {\rm{V_2}}$, then $\delta_{n}^k(\tilde{P})({\rm{V_1}}, {\rm{V_2}}) := n-2+\max\{2, d_{n}^2(\tilde{N})\}+\min\{\delta_{n}^k(\tilde{Q})({\rm{V_1}}, {\rm{V_2}}),\delta_{n}^k(\tilde{N})({\rm{V_1}}, {\rm{V_2}}) \}$, which can be rewritten as $\delta_{n}^k(\tilde{P})({\rm{V_1}}, {\rm{V_2}}) := n-2+\left\lceil{\frac{2 \cdot n}{3}}\right\rceil+\min\{\delta_{n}^k(\tilde{Q})({\rm{V_1}}, {\rm{V_2}}),\delta_{n}^k(\tilde{N})({\rm{V_1}}, {\rm{V_2}}) \}$ \cite{15, 20} so that (as long as $n > 4$) a sufficient but not necessary condition for having $\delta_{n}^k(\tilde{P})({\rm{V_1}}, {\rm{V_2}}) = n-2+\left\lceil{\frac{2 \cdot n}{3}}\right\rceil+\delta_{n}^k(\tilde{Q})({\rm{V_1}}, {\rm{V_2}})$ is trivially given by $\sum_{j=1}^{k}\frac{\left|(x_j+c_j)-x_j \right|}{3} \geq k$, and thus $\sum_{j=1}^{k} \left| c_j \right| \geq 3 \cdot k$ is enough to let us forget the promotion to a knight issue for the given case.

\sloppy Basically (assuming, as usual, $k \in \mathbb{N}-\{0,1\}$ and $n \in \mathbb{N}-\{0,1,2,3,4\}$), we could imagine the $k$-pawn as if it were a human-like piece living in a station (belonging to the $k$-lattice network $C(n,k)$) located at ${\rm{V_1}} \equiv (x_1,x_2,x_3,\ldots,x_k)$ and needed to go in $(x_1,n-1,x_3,\ldots,x_k)$, where he has to alternatively buy a ticket for a $k$-knight-bus or a $k$-queen-train (we assume that both of them travel at the same constant speed) to his final destination at ${\rm{V_2}} \equiv (x_1+c_1,x_2+c_2,x_3+c_3,\ldots,x_k+c_k)$, where ${\rm{V_2}} \in C(n,k) : {\rm{V_2}} \not\equiv {\rm{V_1}}$, making a stop at ${\rm{V_1}}$ first.

Since he cannot take any risk of missing the $k$-knight-bus/$k$-queen-train, let us fairly assume that our $k$-pawn needs to set in advance a comfortable amount of time units to be sure to complete his journey to buy the ticket first (so he sets $n-2$ time units to be certain to reach $(x_1,n-1,x_3,\ldots,x_k)$ in time), at that point he will be faced with two mutually exclusive options: taking the $k$-knight-bus or the $k$-queen-train, and both of them will make a stop at his home station after $\left\lceil{\frac{2 \cdot n}{3}}\right\rceil$ time units by the moment he will have bought the ticket.

The human-like $k$-pawn wishes to reach $(x_1+c_1,x_2+c_2,x_3+c_3,\ldots,x_k+c_k)$ as soon as possible, but he already knows that he cannot change the fact that after $n-2+\left\lceil{\frac{2 \cdot n}{3}}\right\rceil$ time units from the start he will be at his home station, in $(x_1,x_2,x_3,\ldots,x_k)$, sitting on a $k$-knight-bus or a $k$-queen-train.

Thus, he will choose to buy the $k$-knight-bus or $k$-queen-train ticket, depending on ${\rm{V_2}}$.

Consequently, the $k$-pawn diameter problem is to find the minimum amount of time units that our human-like $k$-pawn cannot avoid spending to be able to reach any place in the network $C(n,k)$ in the worst-case scenario (for both ${\rm{V_1}}$ and ${\rm{V_2}}$), by considering his one-time choice at $(x_1,n-1,x_3,\ldots,x_k)$.

Here, we may also consider the chance that a pawn is promoted to a knight instead of a queen since it is allowed by the FIDE Handbook \cite{1} and real chess players decide to promote to a knight sometimes (we can safely disregard the hypothetical promotion to a bishop or a rook by observing that stalemates are not an option in our environment).
As an example, let $n = 4 \wedge k=2$ be given, so the $2$-pawn can start from any vertex $(x_1,x_2) \in C(4,2)$, reach $(x_1,3)$ (if $x_2=3$, it does not matter), promote to a $2$-knight, and then come back to $(x_1,x_2)$ in no more than $4-2+5$ moves (since the maximum knight distance in $C(4,2)$ is $5$, and it is equal to $\delta_{4}^2(\tilde{N})((0,0),(0,3))$). Finally, the promoted piece moves to any other vertex in the set $\{\{0,1,2,3\}\times\{0,1,2,3\}\}-\{(x_1,x_2)\}$ as a $2$-knight or as a $2$-queen, depending on which of the two distances (i.e., the number of required moves as a chess knight or queen) is the shorter.

\begin{figure}[H]
\begin{center}
\includegraphics[width=\linewidth]{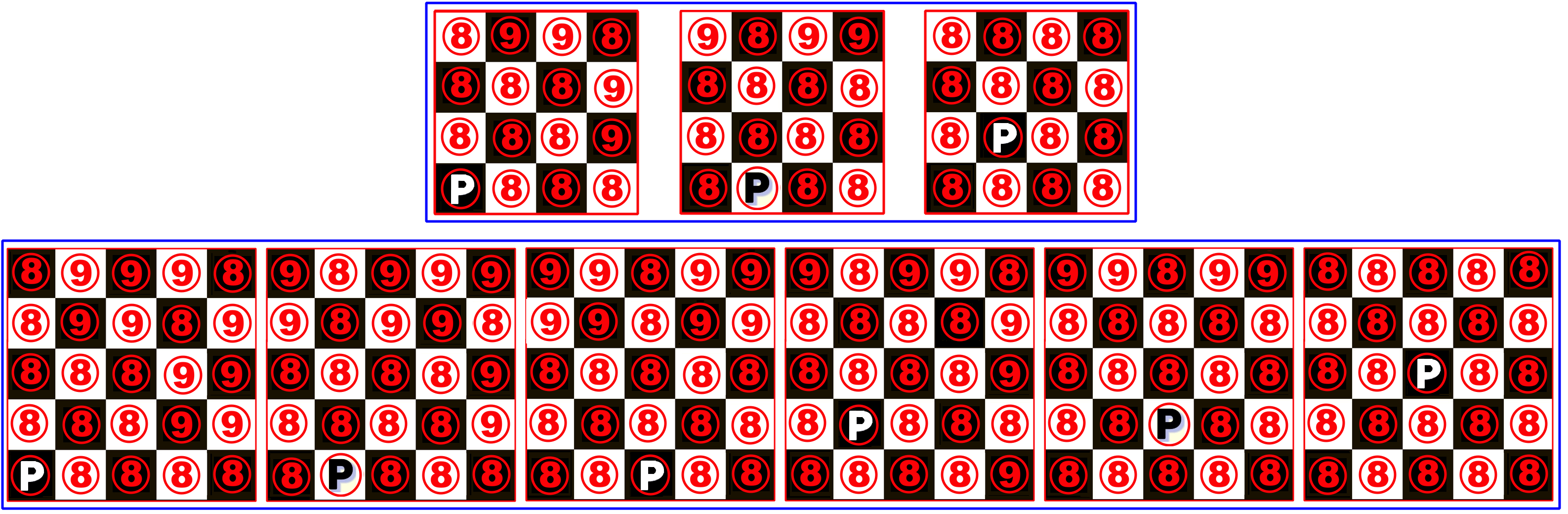}
\end{center}
\caption{Showing by brute force that $r_{4}^2(\tilde{P})=r_{5}^2(\tilde{P})=8$, whereas $d_{4}^2(\tilde{P})=d_{5}^2(\tilde{P})=9$.}
\label{fig:Figure_40}
\end{figure}

\newpage
As a consequence, we get the counterintuitive result $r_{4}^2(\tilde{P})=4-2+5+1=8$, $d_{4}^2(\tilde{P})=4-2+5+2=9$, $r_{5}^2(\tilde{P})=5-2+4+1=8$, and $d_{5}^2(\tilde{P})=5-2+4+2=9$ (see Figure \ref{fig:Figure_40}).

Now, it is clear that the $k$-pawn $\tilde{P}$ also induces a metric on $C(n,k)$; it produces a graph whose circle radius, for any $n > 4$, is smaller than or equal to $n-2+\left\lceil{\frac{2 \cdot n}{3}}\right\rceil+\min\{r_{n}^k(\tilde{Q}),r_{n}^k(\tilde{N}) \}$, and whose diameter, $d_{n}^k(\tilde{P})$, cannot exceed $n-2+\left\lceil{\frac{2 \cdot n}{3}}\right\rceil+\min\{d_{n}^k(\tilde{Q}),d_{n}^k(\tilde{N}) \}$.
Figure \ref{fig:Figure_15} shows that $d_{3}^2(\tilde{P})=3-2+2+1=4$, which is clearly smaller than $3-2+d_{3}^2(\tilde{Q})+\min\{d_{3}^2(\tilde{Q}), d_{3}^2(\tilde{N})\}=3-2+2+\min\{2, +\infty\}=3+2=5$.

Although the $k$-pawn we have just defined is not very similar to the piece shown in Figure~\ref{fig:Figure_8}, we are free to introduce another version of a generalized $k$-dimensional pawn that does not produce a metric (since it can only go forward), but whose moving rule is fully consistent with what we have previously said about the $k$-rook pattern and the light to dark vertex moves requirement (i.e., the parity argument \cite{25}). Hence, let us assume that at the beginning the pawn does not lie on the last rank and then we can state that $S_1(\tilde{P};n,k;x_1,x_2,x_3,\ldots,x_k) = S_1(\overline{P};n,k;x_1,x_2,x_3,\ldots,x_k) = \{(x_1,x_2,x_3,\ldots,x_k)\} \cup \{(x_1, x_2+c, x_3, \ldots, x_k) : c \in \{1,2\} \wedge \left(x_1, x_2+c, x_3, \ldots, x_k \in \{0,1,\ldots,n-1 \}\right)\}$ would be the strictest $k$-dimensional generalization of the standard FIDE pawn to a coherent mathematical environment (Figure \ref{fig:Figure_17}) that does not require any metric space induced by the pawn itself.

At this point, we observe that $\tilde{B}$, $\tilde{R}$, $\tilde{Q}$, $\tilde{P}$ are well-defined and we only need to specify that $\tilde{K}$ and $\tilde{N}$ have been already fully described in Section \ref{sec:Intr} (by Figure \ref{fig:Figure_6} and Figure \ref{fig:Figure_7}, respectively) so that $\tilde{N}:=N$ and $\tilde{K}:=K$.

Lastly, we are free to complete the definition of $\tilde{B}$ as already done in Section \ref{sec:Intr} (i.e., by admitting that a light-square bishop can reach any given dark vertex of $C(n,k)$ at distance $2^k$ and also a bishop placed on any dark vertex of $C(n,k)$ can do the opposite at the same distance) if we wish to get a metric space for the bishop too. If this is the case, then $d_{n}^k(\tilde{B}) \leq 2^k$ clearly holds for any $n,k \in \mathbb{N}-\{0,1\}$.

\subsection{Radius and diameter}
\label{sec:sub4.2}

We are finally ready to bound the values of the radius and diameter for all the six generalized chess pieces, as they have been defined in the previous subsection.

Let us start from the king and assume that $n,k \in \mathbb{N}-\{0,1\}$. Since the definition of $\tilde{K}$ matches the definition of $K$ that has been stated in Section \ref{sec:Intr} (i.e., $S_1(K;n,k;x_1,x_2,\dots,x_k) = \{(x_1+c_1,x_2+c_2,\dots,x_k+c_k) : c_1, c_2, \dots, c_k \in \{-1,0,1\} \hspace{0.5mm} \wedge \hspace{0.5mm} (x_j+c_j) \in \{0,1,\dots,n-1\} \hspace{1mm} \linebreak \forall j \in \{1,2,\dots,k \}\}$, see Figure \ref{fig:Figure_6}), it trivially follows that $d_{n}^k(\tilde{K})=n-1$.
Then, we observe that the $k$-king graph radius corresponds to the maximum number of moves needed by the considered piece to reach the farthest corner of the board from its starting position. Thus, it is optimal to place the $k$-king as close to the center as we can (i.e., if $n$ is even, then it is sufficient to choose the vertex $\left(\frac{n}{2}, \frac{n}{2}, \dots, \frac{n}{2}\right)$ as our starting point, otherwise we will simply start from $\left(\frac{n-1}{2}, \frac{n-1}{2}, \dots, \frac{n-1}{2}\right)$).
Consequently, we have that $r_{n}^k(\tilde{K})=\ceil{\frac{n-1}{2}}$ (in general, given any $k$-chess piece, $X$, we have that $2 \cdot r_{n}^k(\tilde{X}) \geq d_{n}^k(\tilde{X})$ follows from the constraint mentioned in Section \ref{sec:Intr} \cite{11}).

About the $k$-bishop, we simply observe that $d_{n}^k(\tilde{B})=d_{n}^k(B)=r_{n}^k(\tilde{B})=r_{n}^k(B)=2^k$ \linebreak ($n,k \in \mathbb{N}-\{0,1\}$) follows by construction of the bishop metric space (as stated at the end of Subsection \ref{sec:sub4.1}).

The $k$-rook is also a very simple piece to manage (i.e., $\tilde{R} := \overline{R}$) since it changes only one of its $k$ Cartesian coordinates at a time, and this trivially implies that $d_{n}^k(\tilde{R})=r_{n}^k(\tilde{R})=k$, $n,k \in \mathbb{N}-\{0,1\}$.

The diameter of the $k$-knight graph is known as long as $k < 4$.

\begin{figure}[H]
\begin{center}
\includegraphics[width=\linewidth]{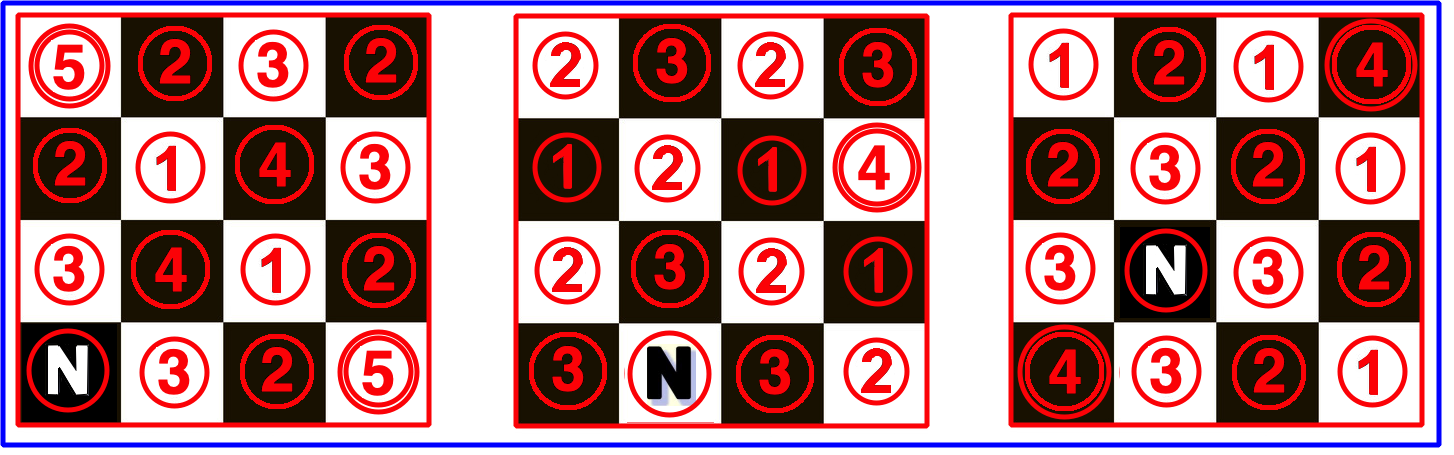}
\end{center}
\caption{Graphical proof that $d_{4}^2(\tilde{N})=5$ and $r_{4}^2(\tilde{N})=4$, since $e_{4}^2(\tilde{N}; 0,0)=5$, $e_{4}^2(\tilde{N}; 1,0)=4$, and $e_{4}^2(\tilde{N}; 1,1)=4$.}
\label{fig:Figure_20}
\end{figure}

In particular, if $k=2$, then Barker's formula (\ref{eq7})
\begin{equation}\label{eq7}
d_{n}^2(\tilde{N}) = \begin{cases}
4, & \text{iff $n \in \{5,6\}$;}\\
5, & \text{iff $n \in \{4,7\}$;}\\
6, & \text{if  $n=8$;} \\
d_{n-1}^2(\tilde{N})+d_{n-3}^2(\tilde{N})-d_{n-4}^2(\tilde{N}), & \text{if  $ n \in \mathbb{N}-\{0,1,2,3,4,5,6,7,8\}$.}
\end{cases}\textnormal{,}
\end{equation}
stated in \cite{15}, gives us the value of $d_{n}^2(\tilde{N})$ for any sufficiently large integer $n$ (for this purpose, Figure \ref{fig:Figure_20} shows that $d_{4}^2(\tilde{N})=5$ is strictly greater than $r_{4}^2(\tilde{N})$), while we have already seen that $\left(k=3 \wedge n \in \mathbb{N}-\{0,1,2,3\}\right) \Rightarrow d_{n}^3(\tilde{N})=n$ (see \cite{16}).

As we can see from Figure \ref{fig:Figure_20}, the knight distance on the graph is a little peculiar since $\delta_{4}^2(\tilde{N})((0,0), (1,1)) = \delta_{4}^2(\tilde{N})((0,0), (2,2))=2 \cdot \delta_{4}^2(\tilde{N})((0,0), (3,3))$, which is quite counterintuitive.

In general, by considering the Manhattan distance between the vertex $(0,0,\ldots,0)$ and the opposite vertex of $C(n,k)$ (i.e., the point $(n-1,n-1,\ldots,n-1)$), we can easily show that the lower bound given by Equation (\ref{eq8})

\begin{equation}\label{eq8}
k \cdot (n-1) \leq d_{n}^k(\tilde{N})\cdot (2+1) \Rightarrow d_{n}^k(\tilde{N}) \geq \left\lceil\dfrac{k \cdot (n-1)}{3}\right\rceil
\end{equation}
holds for any $k \geq 2 \wedge n \geq 4$.

By using the same argument as above (the one we have discussed for the $r_{n}^k(\tilde{K})$ case), we can also find a general lower bound for $r_{n}^k(\tilde{N})$.
Assuming that $n$ is even for the sake of simplicity, let us set the vertex $\left(\frac{n}{2}, \frac{n}{2}, \dots, \frac{n}{2}\right)$ as our starting point (if $n$ is odd, then we will let the $k$-knight start from the vertex $\left(\frac{n-1}{2}, \frac{n-1}{2}, \dots, \frac{n-1}{2}\right)$) and so we have that the Euclidean distance between $\left(\frac{n}{2}, \frac{n}{2}, \dots, \frac{n}{2}\right)$ and $(0,0,\ldots,0)$ is equal to $\sqrt{k} \cdot \frac{n}{2}$ (for any odd value of $n$, the distance is $\sqrt{k} \cdot \frac{n-1}{2}$).
Now, we already know how to improve the trivial lower bound arising from the Euclidean distance between the pair of vertices mentioned above by considering their taxicab distance \cite{3}.
Accordingly, we observe that, at any move (by definition), the $k$-knight adds/subtracts $2$ to one of its initial $k$ Cartesian coordinates and adds/subtracts $1$ to another one so that the following inequality must be satisfied for any $n > 3$

\begin{equation}\label{eq9}
\frac{k \cdot (n-1)}{2} \leq r_{n}^k(\tilde{N}) \cdot (2+1).
\end{equation}

Hence,
\begin{equation}\label{eq19}
r_{n}^k(\tilde{N}) \geq \left\lceil\dfrac{k \cdot (n-1)}{6}\right\rceil.
\end{equation}

On January 2023, the unconstrained (i.e., disregarding the $0 \leq x_1, x_2, \ldots, x_k \leq n-1$ condition) value of $e_{n}^k(\tilde{N}; 0, 0, \ldots, 0)$ has been calculated by the \textit{Mathematics Stack Exchange} user Thomas Andrews as an answer to the thread entitled ``Queen's graph diameter and knight's graph diameter for an $n \times n \times \cdots \times n$ chessboard''.

The aforementioned Andrews' answer implies that
\begin{equation}\label{eq10}
e_{n}^k(\tilde{N}; 0, 0, \ldots, 0) \geq \begin{cases}2+\left\lfloor \frac {k \cdot (n-1)}3\right\rfloor \quad \textnormal{iff} \quad k \cdot (n-1)\equiv 2\pmod 3 \vspace{3mm}; \\
1+\left\lfloor\frac{k \cdot (n-1)}3\right\rfloor \quad \textnormal{iff} \quad k \cdot (n-1)\equiv \{0,1\} \pmod 3.
\end{cases}\textnormal{.}
\end{equation}

Since $e_{n}^k(\tilde{N}; 0, 0, \ldots, 0) \leq d_{n}^k(\tilde{N})$ by definition, the tighter lower bound
\begin{equation}\label{eq11}
d_{n}^k(\tilde{N}) \geq \begin{cases}2+\left\lfloor \frac {k \cdot (n-1)}{3}\right\rfloor \quad \textnormal{iff} \quad k \cdot (n-1)\equiv 2\pmod 3 \vspace{3mm}; \\
1+\left\lfloor\frac{k \cdot (n-1)}{3}\right\rfloor \quad \textnormal{iff} \quad k \cdot (n-1)\equiv \{0,1\} \pmod 3.
\end{cases}
\end{equation}
follows from Equation (\ref{eq10}), for every $n \geq 4$.

The next big challenge is to determine the values of $d_{n}^k(\tilde{Q})$ and $r_{n}^k(\tilde{Q})$ for any $k,n \in \mathbb{N}-\{0,1\}$.

\begin{figure}[H]
\begin{center}
\includegraphics[width=\linewidth]{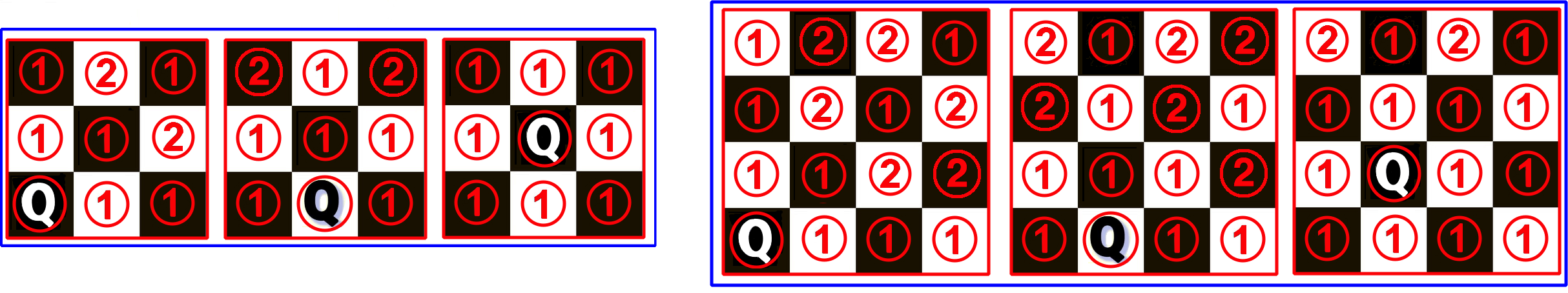}
\end{center}
\caption{Graphical proof that $r_{n}^2(\tilde{Q})=2$ if and only if $n \in \mathbb{N}-\{0,1,2,3\}$, whereas $d_{n}^2(\tilde{Q})=2$ if and only if $n \in \mathbb{N}-\{0,1,2\}$.}
\label{fig:Figure_21}
\end{figure}

If $k=2$, then $\tilde{Q}$ is indistinguishable from $Q$, so $r_{n}^2(\tilde{Q})=2 \hspace{2mm} \forall n \in \mathbb{N}-\{0,1,2,3\}$ and $d_{n}^2(\tilde{Q})=2 \hspace{2mm} \forall n \in \mathbb{N}-\{0,1,2\}$ (see Figure \ref{fig:Figure_21}).

In general, we know that $d_{n}^k(\tilde{R}) \leq k$ holds by definition, and thus $d_{n}^k(\tilde{Q}) \leq k$ by construction. On the other hand, we observe that $r_{n}^k(\tilde{Q}) \geq 2$ for any $n \geq 2 \wedge k \geq 3$ since we need to spend at least two $k$-queen moves to go from $(x_1, x_2, \ldots, x_k)$ to $(x_1+1, x_2+1, \ldots, x_k+1)$, even if we set $k=3$ and we are free to choose any $k$-tuple $(x_1, x_2, \ldots, x_k) : x_1, x_2, \ldots, x_k \in \{0,1,\ldots,n-2\}$.

Hence, $2 \leq r_{n}^k(\tilde{Q}) \leq k$ and $r_{n}^k(\tilde{Q}) \leq d_{n}^k(\tilde{Q}) \leq k$, for any $n >1 \wedge k>2 $.

In this paper, we find $r_{n}^3(\tilde{Q})$ and $d_{n}^3(\tilde{Q})$ by brute force, given the fact that $k:=2 \Rightarrow \tilde{Q}:=Q \Rightarrow \left(r_{n}^k(\tilde{Q})=r_{n}^k(Q)\right) \wedge \left(d_{n}^k(\tilde{Q})=d_{n}^k(Q)\right)$ (see Section \ref{sec:2} and Figure \ref{fig:Figure_21} below).

Thus, let us assume that $k=3$. Our goal is to find which is the minimum value of $n$ such that $r_{n}^3(\tilde{Q})=3$ since $r_{n}^k(\tilde{Q}) \leq d_{n}^k(\tilde{Q}) \leq k$. Then, we can check all those configurations one by one in order to find the exact values of $d_{n}^3(\tilde{Q})$ and $r_{n}^3(\tilde{Q})$.

\begin{figure}[H]
\begin{center}
\includegraphics[width=\linewidth]{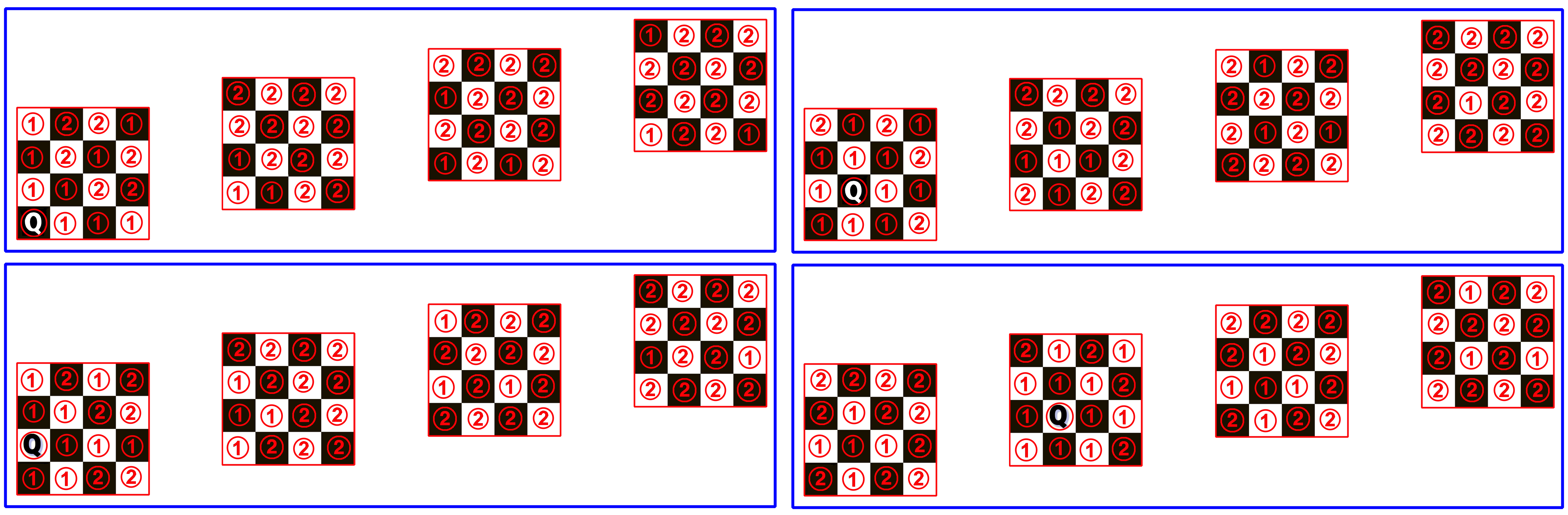}
\end{center}
\caption{Graphical proof that $d_{n}^3(\tilde{Q})=2$ as long as $n \in \{2,3,4\}$.}
\label{fig:Figure_23}
\end{figure}

\begin{figure}[H]
\begin{center}
\includegraphics[width=\linewidth]{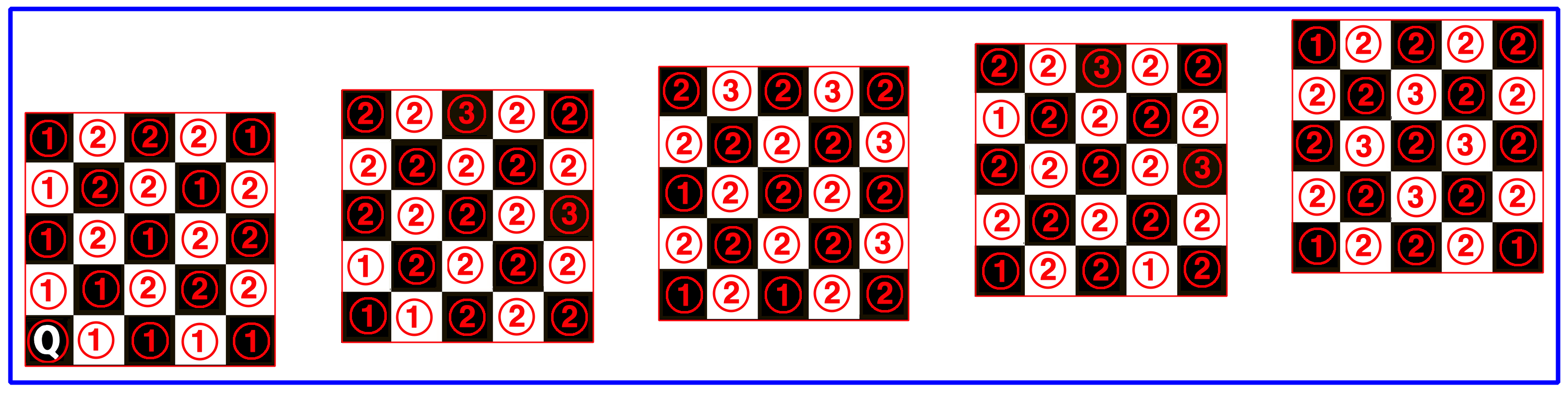}
\end{center}
\caption{Graphical proof that $d_{n}^3(\tilde{Q})=3$ for any $n \in \mathbb{N}-\{0,1,2,3,4\}$.}
\label{fig:Figure_24}
\end{figure}

\begin{figure}[H]
\begin{center}
\includegraphics[width=\linewidth]{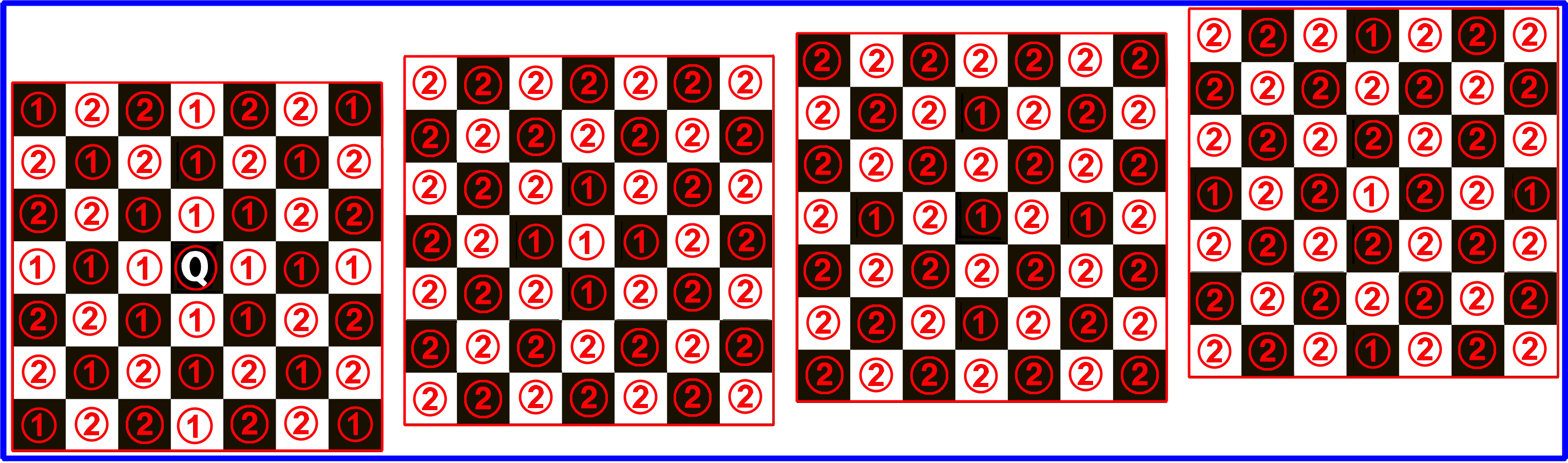}
\end{center}
\caption{Graphical proof that $r_{n}^3(\tilde{Q})=2$ as long as $n \in \{2,3,4,5,6,7\}$ (here we have set $S_0(\tilde{Q};7,3;3,3,3) :=\{(3,3,3)\}$).}
\label{fig:Figure_25}
\end{figure}

\begin{figure}[H]
\begin{center}
\includegraphics[width=\linewidth]{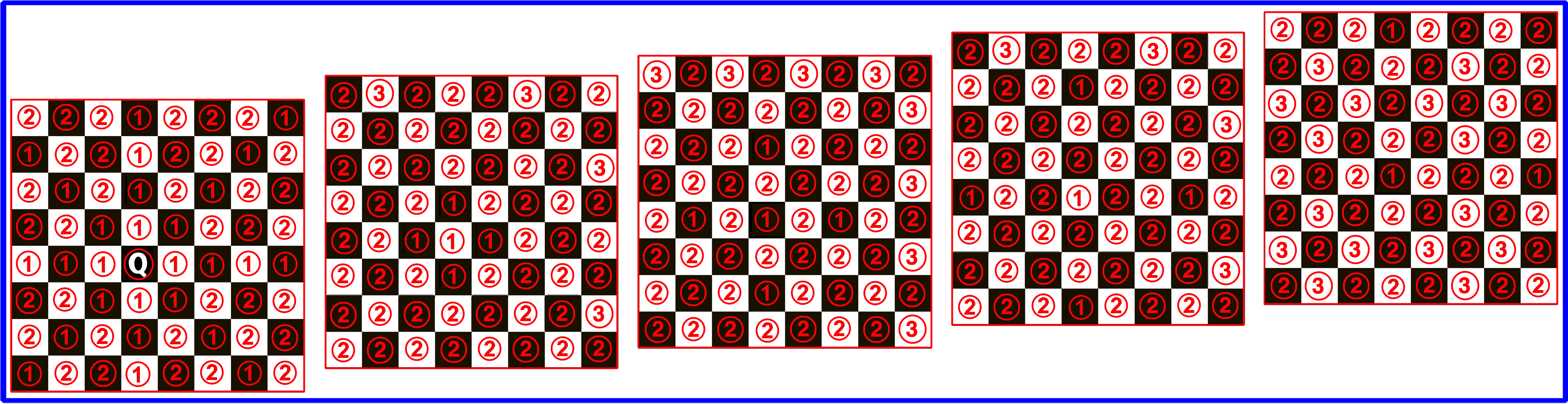}
\end{center}
\caption{Graphical proof that $r_{n}^3(\tilde{Q})=3$ for any $n \in \mathbb{N}-\{0,1,2,3,4,5,6,7\}$ (here we have set $S_0(\tilde{Q};8,3;3,3,3) :=\{(3,3,3)\}$).}
\label{fig:Figure_26}
\end{figure}

As a result (see Figures \ref{fig:Figure_23} to \ref{fig:Figure_26}), we have that 
\begin{equation}\label{eq12}
d_{n}^3(\tilde{Q}) = \begin{cases} 0 \quad \text{iff} \quad n=1; \\
2 \quad \textnormal{iff} \quad n \in \{2,3,4\};\\
3 \quad \textnormal{iff} \quad n \in \mathbb{N}-\{0,1,2,3,4\}.\\
\end{cases}
\end{equation}
and 
\begin{equation}\label{eq13}
r_{n}^3(\tilde{Q}) = \begin{cases} 0 \quad \text{iff} \quad n=1; \\
2 \quad \textnormal{iff} \quad n \in \{2,3,4,5,6,7\};\\
3 \quad \textnormal{iff} \quad n \in \mathbb{N}-\{0,1,2,3,4,5,6,7\}.\\
\end{cases}\textnormal{.}
\end{equation}

Thus, the only remaining $k$-piece to be investigated is $\tilde{P}$.
Since in the previous subsection we have already introduced the tight bounds 
$d_{n}^k(\tilde{P}) \leq n-2+\left\lceil{\frac{2 \cdot n}{3}}\right\rceil+\min\{d_{n}^k(\tilde{Q}),d_{n}^k(\tilde{N}) \}$ and $r_{n}^k(\tilde{P}) \leq n-2+\left\lceil{\frac{2 \cdot n}{3}}\right\rceil+\min\{r_{n}^k(\tilde{Q}),r_{n}^k(\tilde{N}) \}$, we may solve here all the $k=2$ and $k=3$ configurations, leaving the challenge of improving the previous results for any $k \geq 4$ as an open problem for everybody who wishes to have fun with this kind of research.

\begin{figure}[H]
\begin{center}
\includegraphics[width=\linewidth]{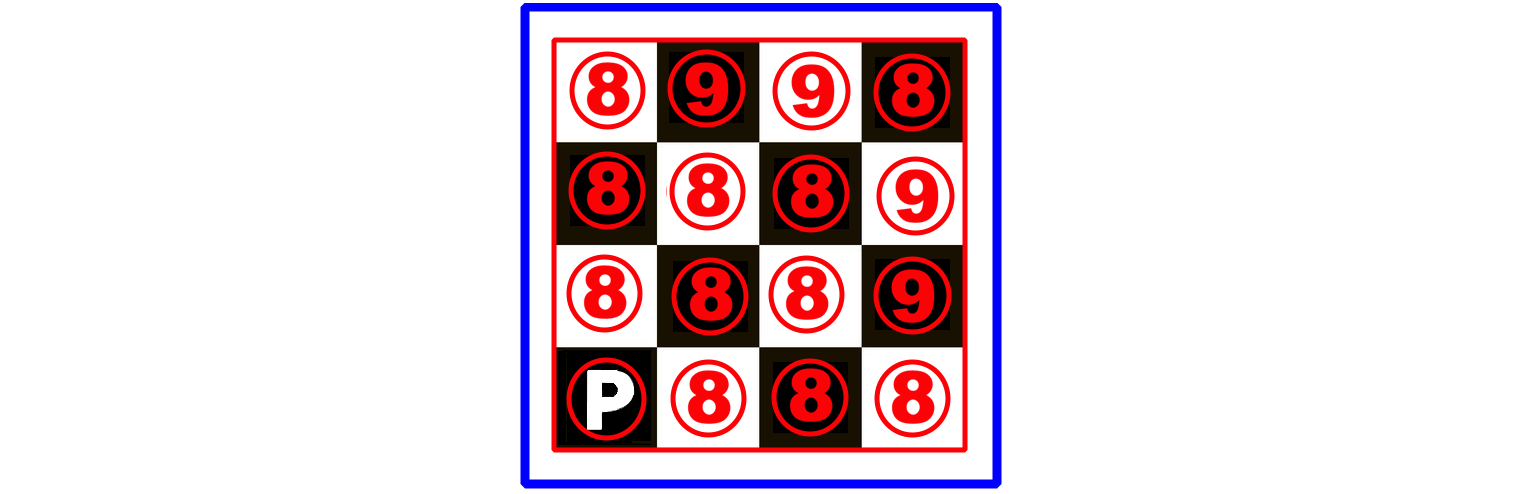}
\end{center}
\caption{Graphical proof that, for any given $n \geq 4$, $d_{n}^2(\tilde{P})=n-2+d_{n}^2(\tilde{N})+2$.}
\label{fig:Figure_42}
\end{figure}

In this regard, let us assume that $n \geq 4$. Then, for any given $\check{k} \in \mathbb{N}-\{0,1\}$, there may exist a minimum value of $n$, $n_{min}:=n_{min}(\check{k})$, such that $r_{n}^{\check{k}}(\tilde{P})=n-2+d_{n}^2(\tilde{N})+\check{k}$ for every $n \geq n_{min}$, and a similar argument should also hold for the smallest $n$ such that $d_{n}^k(\tilde{P})=n-2+d_{n}^2(\tilde{N})+\check{k}$.
More in detail, if $\check{k}=2$, then it is very easy to check that the minimum value of $n$ such that $d_{n}^2(\tilde{P})=n-2+d_{n}^2(\tilde{N})+2$ is $4$ (it follows from the result revealed by Figure~\ref{fig:Figure_15}, given the fact that it is not possible to go from $(0,0)$ to $(3,1)$ with less than two $k$-queen or $k$-knight moves, as shown in Figure~\ref{fig:Figure_42}).

\begin{figure}[H]
\begin{center}
\includegraphics[width=\linewidth]{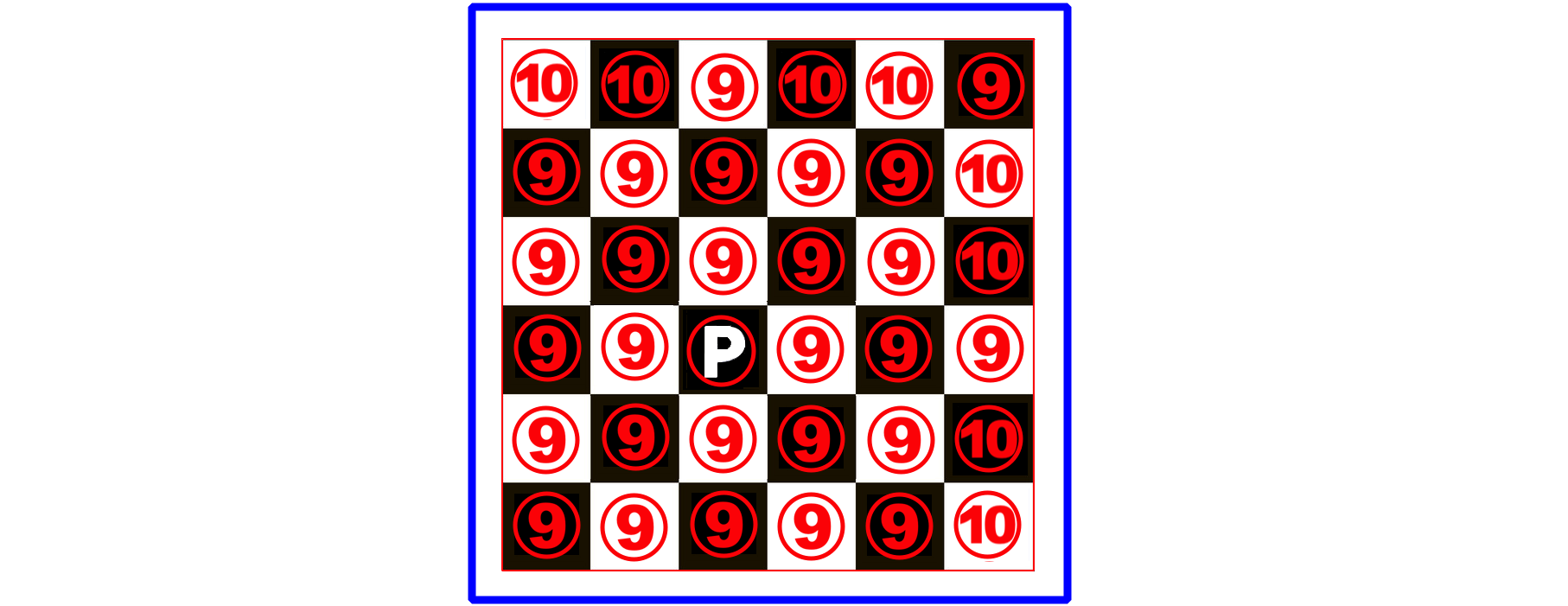}
\end{center}
\caption{One of the six independent configurations that we have checked to prove that $r_{n}^2(\tilde{P})=n-2+d_{n}^2(\tilde{N})+2$ is true for any given $n \geq 6$.}
\label{fig:Figure_43}
\end{figure}

About the minimum $n$ such that $r_{n}^2(\tilde{P})=n-2+d_{n}^2(\tilde{N})+2$, we can see that it is possible to start from the center of the grid (or from $(\frac{n}{2}-1,\frac{n}{2}-1)$ if $n$ is even) and visit with one $k$-queen or $k$-knight move any other point of the set $\{\{0,1,\ldots,n-1\}\times \{0,1,\ldots,n-1\}\}$ as long as $n \leq 5$, whereas if $n=6$ it is simple to check that all the corresponding six independent configurations cannot be covered in just one $k$-queen or $k$-knight move (i.e., if we start from the vertex $(2,2)$, it is trivial to note that we cannot jump to $(5,1)$ with the $2$-knight or reach the aforementioned vertex by moving the $2$-queen only once). This is clearly shown in Figure \ref{fig:Figure_43}.

Therefore, $d_{n}^2(\tilde{P})=n-2+d_{n}^2(\tilde{N})+2$ for any $n \geq 4$ and $r_{n}^2(\tilde{P})=n-2+d_{n}^2(\tilde{N})+2$ if and only if $n \geq 6$.

We can now move on $k=3$. Our strategy is very similar to the $k=2$ case since the values of $n_{min}$, for $d_{n}^3(\tilde{P})=n-2+d_{n}^2(\tilde{N})+3$ and $r_{n}^3(\tilde{P})=n-2+d_{n}^2(\tilde{N})+3$, are small enough to be attacked by brute force.

\begin{figure}[H]
\begin{center}
\includegraphics[width=\linewidth]{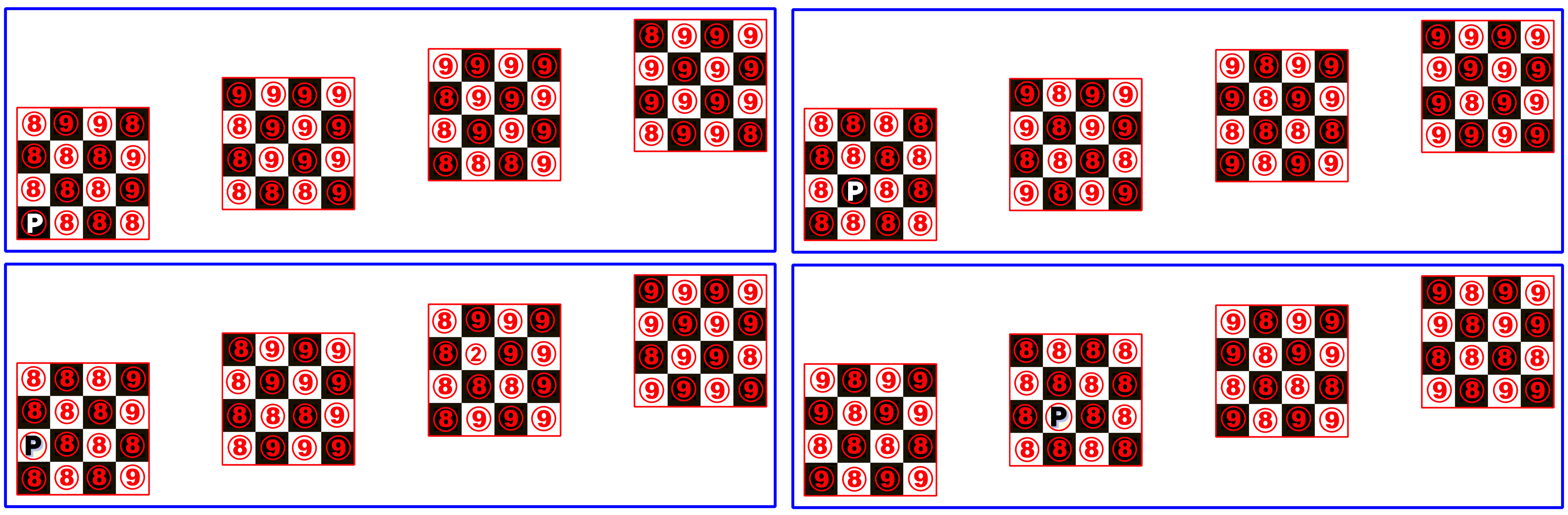}
\end{center}
\caption{Graphical proof that $d_{4}^3(\tilde{P})<4-2+d_{4}^2(\tilde{N})+3$. Thus, the minimum value of $n$ such that $d_{n}^3(\tilde{P})=n-2+d_{n}^2(\tilde{N})+3$ must be greater than or equal to $5$ (see Figure \ref{fig:Figure_45} below in order to verify that $d_{n}^3(\tilde{P})=n-2+d_{n}^2(\tilde{N})+3$ holds for any $n \in \mathbb{N}-\{0,1,2,3,4\}$).}
\label{fig:Figure_44}
\end{figure}

\begin{figure}[H]
\begin{center}
\includegraphics[width=\linewidth]{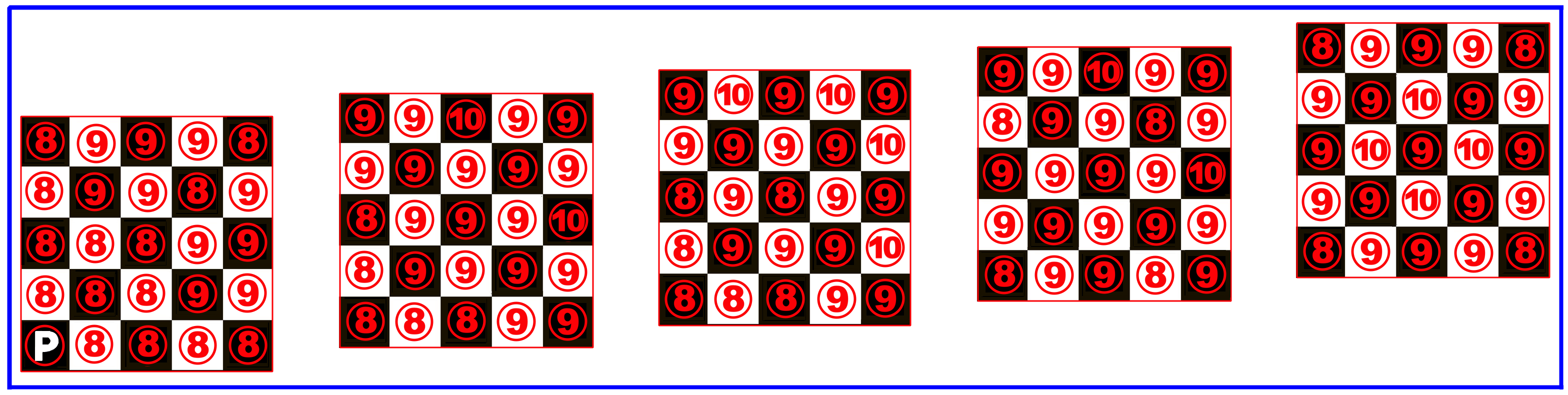}
\end{center}
\caption{Graphical proof that $d_{n}^3(\tilde{P})=n-2+\left\lceil{\frac{2 \cdot n}{3}}\right\rceil+3$ for any $n \in \mathbb{N}-\{0,1,2,3,4\}$.}
\label{fig:Figure_45}
\end{figure}

\begin{figure}[H]
\begin{center}
\includegraphics[width=\linewidth]{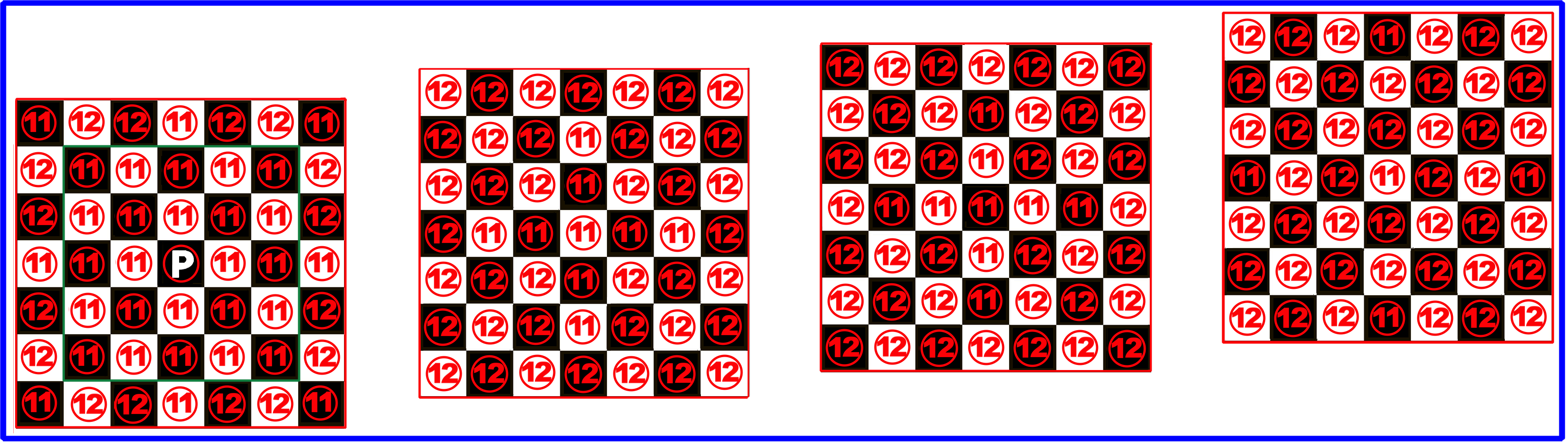}
\end{center}
\caption{Graphical proof that $r_{n}^k(\tilde{P})<n-2+d_{n}^2(\tilde{N})+k$ as long as $k=3$ and $3 < n < 8$ (here we have set $S_0(\tilde{P};7,3;3,3,3) :=\{(3,3,3)\}$).}
\label{fig:Figure_46}
\end{figure}

\begin{figure}[H]
\begin{center}
\includegraphics[width=\linewidth]{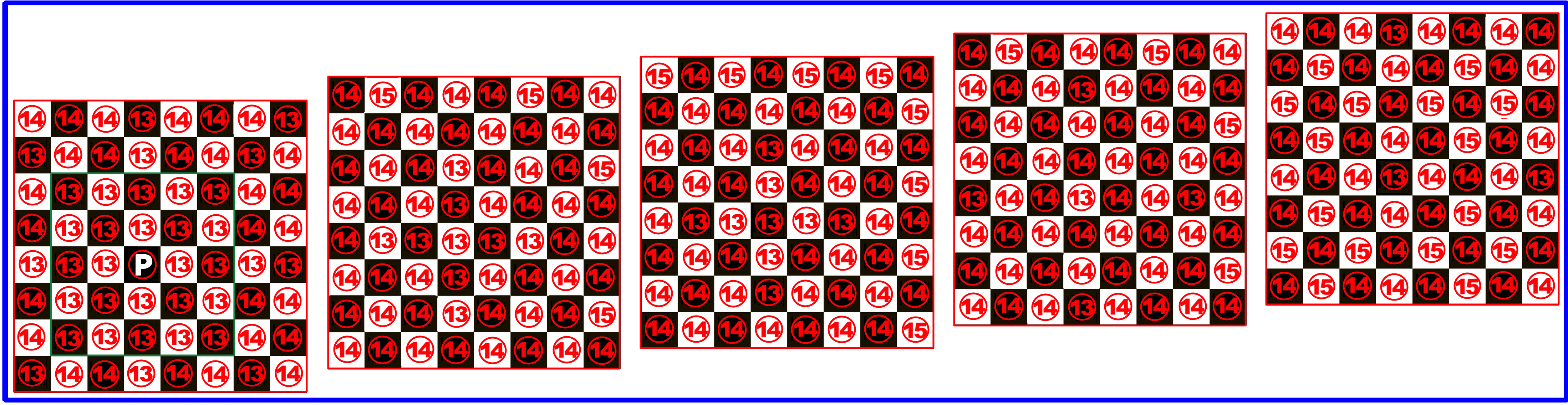}
\end{center}
\caption{One of the $20$ independent configurations that we have checked to prove that $r_{8}^3(\tilde{P})=8-2+6+3=15$ (here we have set $S_0(\tilde{P};8,3;3,3,3) :=\{(3,3,3)\}$) implies $r_{n}^3(\tilde{P})=n-2+\left\lceil{\frac{2 \cdot n}{3}}\right\rceil + 3$ for any $n \in \mathbb{N}-\{0,1,2,3,4,5,6,7\}$.}
\label{fig:Figure_47}
\end{figure}

Figures \ref{fig:Figure_44} to \ref{fig:Figure_47} show that $d_{n}^3(\tilde{P})=n-2+d_{n}^2(\tilde{N})+3$ if and only if $n \geq 5$ (which implies that $d_{n}^3(\tilde{P})=n-2+\left\lceil{\frac{2 \cdot n}{3}}\right\rceil+3$ for any $n > 4$), whereas $r_{n}^3(\tilde{P})=n-2+d_{n}^2(\tilde{N})+3$ for any $n \geq 8$.

In detail, we have that $\delta_{5}^3(\tilde{N})((0,0,0), (4,2,1))=3 \Rightarrow \delta_{5}^3(\tilde{P})((0,0,0), (4,2,1)) = 5-2+\max\{d_{5}^2(\tilde{Q}), d_{5}^2(\tilde{N}) \}+\min\{\delta_{5}^3(\tilde{Q})((0,0,0), (4,2,1)),\delta_{5}^3(\tilde{N})((0,0,0), (4,2,1)) \}=5-2+4+3=10$ and $\delta_{8}^3(\tilde{N})((3,3,3), (7,0,2))=3 \Rightarrow \delta_{8}^3(\tilde{P})((3,3,3), (7,0,2)) = 8-2+\max\{d_{8}^2(\tilde{Q}), d_{8}^2(\tilde{N}) \}+\min\{\delta_{8}^3(\tilde{Q})((3,3,3), (7,0,2)),\delta_{8}^3(\tilde{N})((3,3,3), (7,0,2)) \}=8-2+6+3=15$ (let us omit to write here the remaining $19$ configurations for $k=3 \wedge n=8$ that we have rigorously checked in order to state the above results).


\section{Conclusion} \label{sec:5}

We have defined a coherent set of metric spaces induced by highly reliable generalizations to $k$-dimensions of the six planar FIDE chess pieces. In order to satisfy the triangle inequality, we have stated that the bishop can reach in $2^k$ moves any given vertex with opposite parity from the starting one, while the $k$-pawn will satisfy the condition ${\delta}({\rm{V_2}}, {\rm{V_1}})={\delta}({\rm{V_1}}, {\rm{V_2}}) \hspace{2mm} \forall \hspace{0.5mm}{\rm{V_1}}, {\rm{V_2}}\in C(n,k)$ (i.e., the symmetry property of a metric space) by ideally coming back to ${\rm{V_1}}$ as a promoted $k$-queen or $k$-knight after a sufficiently large number of moves, before finally reaching ${\rm{V_2}}$ traveling through the shortest $k$-queen or $k$-knight path.

Then, by applying retrograde check to the constraints 1, 2, and 3, stated in Subsection \ref{sec:sub4.1}, we have managed to define the moving rules of all the $k$-dimensional extensions of the standard FIDE chess pieces, according to the third article of the FIDE Laws of Chess Handbook \cite{1}.

Subsequently, the move rules describing our $k$-pieces are as follows.
 
$k$-bishop (see Figure \ref{fig:Figure_22}). $S_0(\tilde{B};n,k;x_1,x_2,\dots,x_k):=\{(x_1, x_2, \dots, x_k)\} \rightarrow \linebreak S_1(\tilde{B};n,k;x_1, x_2, \dots, x_k)=\{(x_1, x_2, \dots, x_k)\} \cup \{(x_1+c_1, x_2+c_2, \dots, x_k+c_k) : \linebreak \exists c \in \{1,\dots,n-1\}\quad \forall j \in \{1,2, \dots, k\} \hspace{2mm} c_j \in \{-c, 0, c \} \wedge |\{j : c_j \neq 0 \}|=2 \wedge (x_j+c_j) \in \{0,1,\dots,n-1 \}\}$.

$k$-rook (see Figure \ref{fig:Figure_4}). $S_0(\tilde{R};n,k;x_1,x_2,\dots,x_k):=\{(x_1, x_2, \dots, x_k)\} \rightarrow \linebreak S_1(\tilde{R};n,k;x_1,x_2,\dots,x_k) = \bigcup_{1 \leq j \leq k} \{ (x_1, \ldots, x_{j-1}, c, x_{j+1}, \ldots, x_k) \hspace{2mm} : \hspace{2mm} c \in \{0, 1, \ldots, n-1\} \}$.
 
$k$-queen (see Figure \ref{fig:Figure_30}).
$S_0(\tilde{Q};n,k;x_1,x_2,\dots,x_k):=\{(x_1, x_2, \dots, x_k)\} \rightarrow \linebreak S_1(\tilde{Q};n,k;x_1,x_2,\dots,x_k)=S_1(\tilde{B};n,k;x_1, x_2, \dots, x_k) \cup S_1(\tilde{R};n,k;x_1,x_2,\dots,x_k)$.
 
$k$-king (see Figure \ref{fig:Figure_6}).
$S_0(\tilde{K};n,k;x_1,x_2,\dots,x_k):=\{(x_1, x_2, \dots, x_k)\} \rightarrow \linebreak S_1(K;n,k;x_1,x_2,\dots,x_k) = \{(x_1+c_1,x_2+c_2,\dots,x_k+c_k) : c_1, c_2, \dots, c_k \in \{-1,0,1\} \wedge \linebreak (x_j+c_j) \in \{0,1,\dots,n-1\} \hspace{2mm} \forall j \in \{1,2,\dots,k \}\}$.

$k$-knight (see Figure \ref{fig:Figure_7}).
$S_0(\tilde{N};n,k;x_1,x_2,\dots,x_k):=\{(x_1, x_2, \dots, x_k)\} \rightarrow \linebreak S_1(N;n,k;x_1, x_2, \dots, x_k)=\{(x_1, x_2, \dots, x_k)\} \cup \{(x_1+c_1, x_2+c_2, \dots, x_k+c_k) : \linebreak \forall j \in \{1,2,\dots,k\} \bigl( \sum_{j=1}^{k} |c_j| = 3 \wedge |c_j| \in \{0,1,2\} \wedge \exists! j : 
|c_j|=1 \bigr) \wedge (x_j+c_j) \in \{0,1,\dots,n-1 \}\}$.

About the $k$-pawn, we could simply state that $S_0(\tilde{P};n,k;x_1,x_2,x_3,\ldots,x_k):= \linebreak \{(x_1, x_2, x_3, \ldots, x_k)\} \rightarrow S_1(\tilde{P};n,k;x_1,x_2,x_3,\ldots,x_k) = \{(x_1,x_2,x_3,\ldots,x_k)\} \cup \linebreak \{(x_1, x_2+c, x_3, \ldots, x_k) : c \in \{1,2\} \wedge \left(x_1, x_2+c, x_3, \ldots, x_k \in \{0,1,\ldots,n-1 \}\right)\}$ (see Figure \ref{fig:Figure_8}), but this assumption would not induce any metrics on $C(n,k)$. Consequently, in Subsection \ref{sec:sub4.1}, for any $n \geq 5$, $k \geq 2$, ${\rm{V_1}} \in C(n,k) : {\rm{V_1}} \equiv (x_1, x_2, \ldots, x_k)$, and ${\rm{V_2}} \in C(n,k) : {\rm{V_2}} \equiv (x_1+c_1,x_2+c_2,\ldots,x_k+c_k) \ne {\rm{V_1}}$, we have finally defined the $k$-pawn distance between ${\rm{V_1}}$ and ${\rm{V_2}}$ as $\delta_{n}^k(\tilde{P})({\rm{V_1}}, {\rm{V_2}}) := n-2+\left\lceil{\frac{2 \cdot n}{3}}\right\rceil+\min\{\delta_{n}^k(\tilde{Q})({\rm{V_1}}, {\rm{V_2}}),\delta_{n}^k(\tilde{N})({\rm{V_1}}, {\rm{V_2}}) \}$.

Then, for each of the graphs induced on $C(n,k):=\{\{0,1,\dots,n-1\}\times \{0,1,\dots,n-1\} \times \cdots \times \{0,1,\dots,n-1\}\} \subseteq \mathbb{Z}^k$ by the $k$-pieces above, we have bounded the values of the radius and diameter, fully solving the problem for the $k$-bishop (trivially), $k$-rook, $k$-king, $2$-queen and $3$-queen, and $2$-pawn and $3$-pawn. We also give the diameters of the $2$-knight and $3$-knight graphs (see Subsection \ref{sec:sub4.2}).
Moreover (assuming $k \geq 2$ and $n \geq 5$), we know that $r_{n}^k(\tilde{P}) \leq n-2+\left\lceil{\frac{2 \cdot n}{3}}\right\rceil+r_{n}^k(\tilde{Q})$ and $d_{n}^k(\tilde{P}) \leq n-2+\left\lceil{\frac{2 \cdot n}{3}}\right\rceil+d_{n}^k(\tilde{Q})$.

The problem is open for the remaining cases, keeping in mind that $k \geq d_{n}^k(\overline{R}) \geq d_{n}^k(\tilde{Q}) \geq d_{n}^k(Q)$ holds by definition.

\makeatletter
\renewcommand{\@biblabel}[1]{[#1]\hfill}
\makeatother

\end{document}